\documentclass[a4paper, 11pt]{article}
\usepackage{amsmath,amsfonts,latexsym,theorem}
\usepackage{graphicx}
\nonstopmode
\graphicspath{{figures/}}

\numberwithin{equation}{section}

\newtheorem{Theorem}{Theorem}[section]
\newtheorem{Definition}{Definition}[section]

\newtheorem{Lemma}{Lemma}[section]

\newenvironment{Proofc}[1]{\smallskip\par\noindent\textsc{#1}\quad}%
  {\hfill$\Box$\bigskip\par}
\newenvironment{Proof}{\begin{Proofc}{Proof}}{\end{Proofc}}

\newtheorem{Remark}{Remark}[section]

\def\a{\alpha}
\def\b{\beta}
\def\d{\delta}

\def\g{\gamma}
\def\G{\Gamma}
\def\l{\lambda}

\def\o{\omega}

\def\pd{\partial}

\newcommand{\ds}{\displaystyle}

\newcommand{\cA}{{\cal A}}
\newcommand{\cB}{{\cal B}}
\newcommand{\cC}{{\cal C}}
\newcommand{\cE}{{\cal E}}
\newcommand{\cF}{{\cal F}}
\newcommand{\cG}{{\cal G}}
\newcommand{\cH}{{\cal H}}
\newcommand{\cI}{{\cal I}}
\newcommand{\cK}{{\cal K}}

\newcommand{\cR}{{\cal R}}
\newcommand{\cS}{{\cal S}}
\newcommand{\cT}{{\cal T}}
\newcommand{\cV}{{\cal V}}

\newcommand{\R}{{\mathbb R}}

\def\pd{\partial}

\newcommand{\Inc}{\text{Inc}}
\begin{document}
\title{A numerical method for Mean Field Games on networks}
\author{Simone Cacace\footnotemark[1],\,\, Fabio Camilli\footnotemark[2]\,\, and Claudio Marchi\footnotemark[3] }
\date{}
\vskip -1.5cm

\maketitle

\footnotetext[1]
{
 Dip. di Matematica, ``Sapienza'' Universit{\`a}  di Roma, p.le A. Moro 5, 00185 Roma, Italy,
 ({\tt e-mail: cacace@mat.uniroma1.it})
}
\footnotetext[2]
{
 Dip. di Scienze di Base e Applicate per l'Ingegneria,  ``Sapienza" Universit{\`a}  di Roma,
 via Scarpa 16, 00161 Roma, Italy, ({\tt e-mail: camilli@sbai.uniroma1.it})
}
\footnotetext[3]
{
 Dip. di Ingegneria dell'Informazione, Universit\`a di Padova, via Gradenigo 6/B, 35131 Padova, Italy
 ({\tt e-mail: claudio.marchi@unipd.it})
}
\begin{abstract}
We propose  a  numerical method  for stationary   Mean Field Games defined on a network. In this framework
a correct approximation of the transition conditions at the vertices plays a crucial role. We   prove existence, uniqueness and convergence of the scheme and we also propose a least squares method
for  the solution of the discrete system. Numerical experiments are carried out.
\end{abstract}
 \begin{description}
\item [\textbf{ AMS subject classification}:  ] 91A15, 35R02, 35B30, 49N70, 65M06.
  \item [\textbf{Keywords}: ] networks, mean field games, finite difference schemes, convergence.
 \end{description}
%
%
%
\section{Introduction}\label{intro}
The Mean Field Game (MFG in short) theory   has been   introduced in  \cite{hmc,ll}
to describe the limit behavior of differential games when the number of agents becomes
very large. Models based on this theory can be used to investigate crowd dynamics,
consensus formation and various economical and social problems (as growth theory, environmental policy and formation of volatility in financial markets) in which the strategy
of the single agent determines a collective behavior of the population (see \cite{abllm,do,gll}).\par

From a mathematical  point of view, MFG theory leads to the study of a coupled system of two differential equations: a  Hamilton-Jacobi-Bellman equation and a Fokker-Planck equation, describing respectively the optimal behavior of each single agent and the evolution of the whole population.
There is a rapidly increasing literature concerning both the theoretical aspects and the
applications of MFG (see the review paper \cite{gs}). \par

A crucial point to extend the theory of MFG systems to networks is to find the appropriate
transition conditions at the vertices in order to obtain a well posed mathematical problem,
coherent with the applications. In \cite{ccm}, it was considered a MFG system with quadratic
Hamiltonian which, by an appropriate change of variable, can be transformed into a linear
system of differential equations coupled only   via the initial datum. A general class of stationary MFG systems on networks
is considered in \cite{cm}, where it is proved existence and uniqueness of classical solutions
to the problem
\begin{equation}\label{MFGs}
 \left\{
  \begin{array}{ll}
   -\nu \pd^2 u  +H(x, \pd u)+\l  =V[m],\qquad &x\in \G\\[4pt]
    \nu \pd^2 m +\pd(m\,H_p(x, \pd u))=0,&x\in \G\\[4pt]
    \int_\G m(x)dx=1, \int_\G u(x)dx=0.
  \end{array}
 \right.
\end{equation}
Here the network $\G=(\cV,\cE)$ is a finite collection of points $\cV:=\{v_i\}_{i\in I}$ in $\R^n$
indexed by $I$, connected by continuous, non self-intersecting arcs $\cE:=\{e_j\}_{j\in J}$
indexed by $J$. Moreover, $\nu=(\nu_j)_{j\in J}$ are strictly positive numbers,
the Hamiltonian $H$ is a collection $\{H_j\}_{j\in J}$ where $H_j$ are continuous, convex Hamiltonians
defined on the arcs $e_j$, and $H_p(x,q)$ is the differential of $p\mapsto H(x,p)$ at $p=q$. Let us stress that $H$ may be discontinuous at the vertices. \\
The equations in   \eqref{MFGs}, which have the same interpretation as in the classical MFG theory,
are defined in terms of the coordinate parametrizing the arc $e_j$, $j\in J$, and have to be
complemented with appropriate  conditions at the vertices. At each internal vertex $v_i$ we consider
the transition conditions
\begin{equation} \label{kir}
\begin{split}
   & \sum_{j\in \Inc_i } \nu_j\pd_j u(v_i)=0,\\
    &\sum_{j\in \Inc_i}[ \nu_j\pd_j m(v_i)+H_{j,p }(v_i,\pd_j u) m_j(v_i)]=0,\\
    &u_j(v_i)=u_k(v_i), \,m_j(v_i)=m_k(v_i), \quad j,k\in \Inc_i
\end{split}
\end{equation}
where $\Inc_i$ denotes the set of the edges incident the vertex $v_i$ and $H_{j,p}$ is the derivative
with respect to $p$ of the Hamiltonian  defined on the edge $e_j$.
We mention, see \cite{cm} for more details, that the first condition in \eqref{kir} is the classical
Kirchhoff condition and it prescribes the probability that an agent reaching the vertex $v_i$ enters
in the incident edge $e_j$, $j\in \Inc_i$; the second condition in \eqref{kir} guarantees the mass
conservation at $v_i$ (the sum of the fluxes at $v_i$  is null); the third condition is the continuity
of $u$ and $m$ at $v_i$.
\textbf{We remark that  \eqref{kir} are natural conditions for     $2^{nd}$ order problems on networks. In fact the domain of the Laplace operator on a network  is given by  continuous functions on $\G$ which are $H^2$ on the edges  and which  satisfy the Kirchhoff condition at the vertices
\cite{nic}.  Moreover,  the  transition conditions  are a crucial ingredient for the validity of the  maximum principle on networks.}\par

In this paper we consider the numerical approximation of the problem \eqref{MFGs}-\eqref{kir}
following the approach in \cite{a,acd}, where a finite difference approximation of the MFG system
is studied (see also \cite{cs,g,lst} for different approaches).
Inside the edges we follow the same  approach of \cite{acd} and we discretize the differential
equations in \eqref{MFGs}  by finite differences. The guideline to find the correct approximation
of the transition conditions in \eqref{kir} is to reproduce at a discrete level some fundamental
identities which are obtained in the continuous setting by the weak formulation of the problem
(see f.e. \eqref{fundid}). For this reason the discrete Hamiltonian defined by a monotone approximation
of the Hamilton-Jacobi-Bellman equation is also used in the discretization of the Fokker-Planck equation
and of the corresponding transition condition. By means of the previous identities we prove the
well-posedness of the discrete problem and the convergence to the solution of the system \eqref{MFGs}. \par

\textbf{While there is a large  literature about the  approximation of hyperbolic problems on networks
(see for example \cite{bnc}, \cite{damp}), as far as we know, numerical schemes for second order differential equations on networks with Kirchhoff
conditions have been only considered   in the linear case (see \cite{lz}, \cite{mn})}. Hence the part concerning
the approximation of the Hamilton-Jacobi-Bellman equation on the network is new and of independent interest.\par

The paper is organized as follows. In Section \ref{sec1} we introduce assumptions and notations. Section \ref{sec2}
includes three subsections concerning existence, uniqueness and
convergence. In Section \ref{sec4} we present  a method for the solution of the discrete system and  some numerical examples illustrating the theory.

\section{Notations and preliminary definitions}\label{sec1}
\begin{figure}[h!]
 \begin{center}
 \includegraphics[width=.3\textwidth]{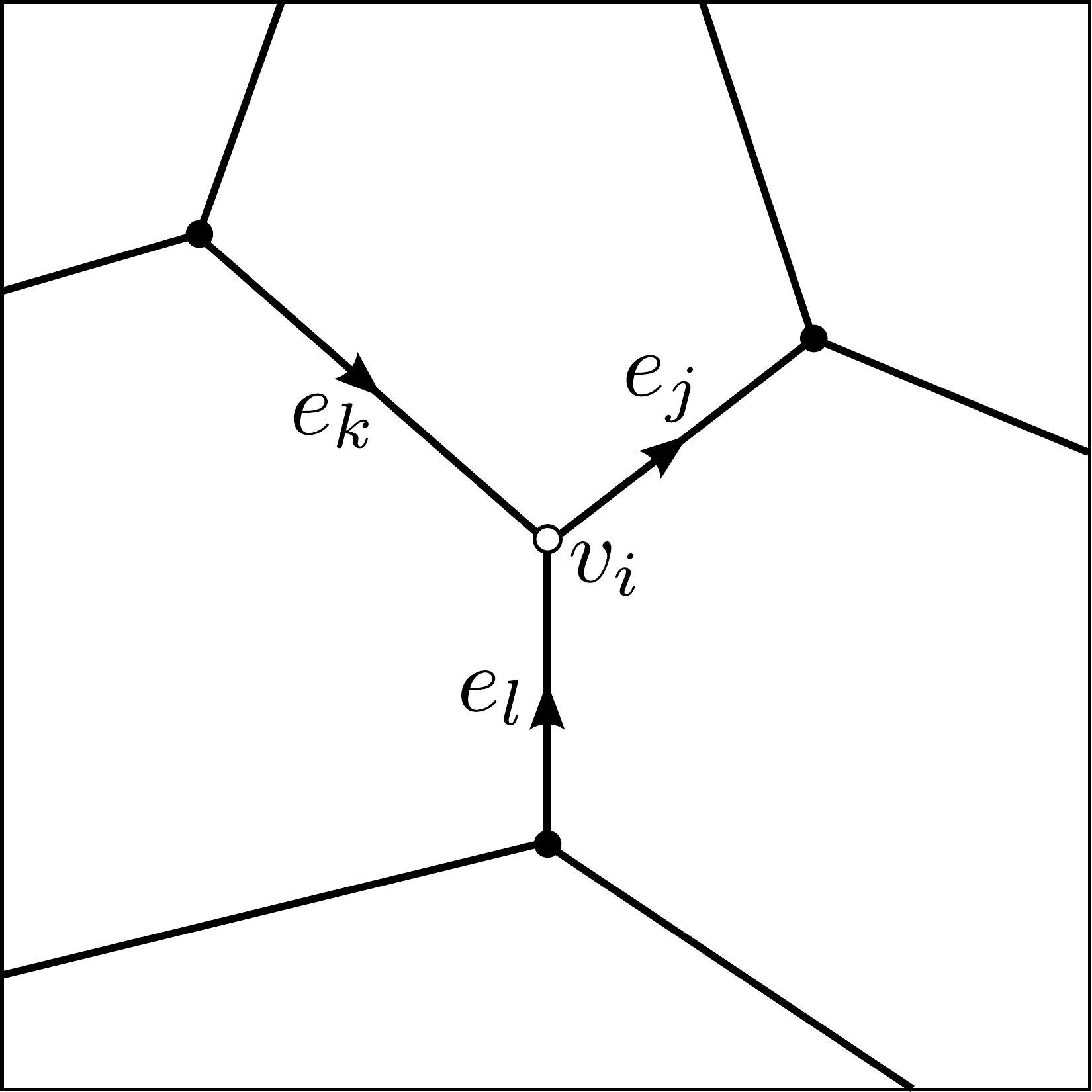}
 \end{center}
 \caption{incident edges to the vertex $v_i$: $\Inc_i^{+}=\{j\}$, $\Inc_i^{-}=\{k,l\}$.}\label{network-inc}
 \end{figure}

A {\it network} is a couple $({\cal V},{\cal E})$ given by a finite collection of vertices ${\cal V}:=\{v_i\}_{i\in I}$ and a finite collection ${\cal E}:=\{e_j\}_{j\in J}$ of continuous non self-intersecting arcs whose endpoints belong to~${\cal V}$.
 We assume that each arc $e_j\in \cE$ is  parametrized by a smooth function $\pi_j:[0,l_j]\to\R^n$, $l_j>0$.
For a function $u:  \G\to\R$  we denote by $u_j:[0,l_j]\to \R$ the restriction of $u$ to $e_j$, i.e.
$u(x)=u_j(y)$ for $x\in e_j$, $y=\pi_j^{-1}(x)$. Given $v_i\in \cV$, we denote by $\pd_j u(v_i)$
the oriented derivative at $v_i$ along the arc $e_j$ defined by
\[
\pd_j u (v_i)=\left\{
               \begin{array}{ll}
                 \displaystyle\lim_{t\to 0^+}(u_j(t)-u_j(0))/t, & \hbox{if $v_i=\pi_j(0)$;} \\
                 \displaystyle\lim_{t\to 0^+}(u_j(l_j-t)-u_j(l_j))/t, & \hbox{if $v_i=\pi_j(l_j)$.}
               \end{array}
             \right.
\]
Given a discretization step $h=\{h_j\}_{j\in J}$, we consider an uniform partition
$y_{j,k}=k h_j$, $k=0,\dots,N^h_j,$  of the interval $[0,l_j]$ which parameterizes the edge $e_j$
(we assume that $N^h_j=l_j/h_j$ is an integer). We obtain a spatial grid on $\G$ by setting
\begin{equation}\label{grid}
  \cG_h=\{x_{j,k}=\pi_j(y_{j,k}),\; j\in J, \; k=0,\dots,N^h_j\}.
\end{equation}
We define
$$
   \Inc_i^{+}=\{j\in \Inc_i:\, v_i=\pi_j(0)\}, \quad \Inc_i^{-}=\{j\in \Inc_i:\, v_i=\pi_j(N^h_jh_j)\},\\
$$ so that
$$ \Inc_i=\Inc_i^{+}\cup \Inc_i^{-}\,,$$
as shown in Figure \ref{network-inc}.
\\
We set
\begin{equation}\label{hi}
|h|=\max_{j\in J}\{h_j\}, \qquad  h_{v_i}=\sum_{j\in \Inc_i}\frac{h_j}{2},\qquad N^h=\#(I)+\sum_{j\in J} (N^h_j-1),
\end{equation}
i.e. $N^h$ is the total number of the points of $\cG_h$, having identified for each $i\in I$
the $\#(\Inc_i)$ grid points corresponding to the same vertex $v_i$.
\textbf{For a grid function $U:\cG_h\to \R$ we denote by  $U_{j,k}$ its value  at the grid point $x_{j,k}$.}
\begin{Definition}
We say that a grid function $U:\cG_h\to \R$ is continuous at $v_i$ if
\[ U_{j,\ell}=U_{k,m}:=U_i\quad  \text{if $v_i=\pi_j(\ell h_j)=\pi_k(m  h_k)$, $ j, k\in \Inc_i$, $\ell\in\{0, N^h_j\}$, and $m\in\{0, N^h_k\}$},\]
 i.e. the value  of $U$ at the   vertex $v_i$ is independent of incident  edge $e_j$, $j\in \Inc_i$.
We say that a a grid function is continuous if it is continuous at $v_i$, for each $i\in I$.
\end{Definition}
 We introduce the  finite difference operators
 \begin{align*}
&(D^+ U)_{j,k} = \frac{ U_{j,k+1}-U_{j,k}   } {h_j},\\
&[D_h U]_{j,k} = \big((D^+ U )_{j,k} , (D^+ U )_{j,k-1}\big)^T,\\
&(D^2_h U)_{j,k}=   \frac { U_{j,k-1}-2U_{j,k}+U_{j,k+1}} {h_j^2}.
\end{align*}
In order to approximate the Hamiltonian $H_j:[0,l_j]\times \R\to\R$, $j\in J$,
we consider a  numerical Hamiltonian $g_j: [0,l_j]\times \R^2\to\R$,\,  $(x,q_1,q_2)\to g_j\left(x,q_1,q_2\right)$
satisfying the following assumptions:\\
\\ ($\mathbf{G1}$)  \emph{monotonicity: $g_j$ is nonincreasing with respect to $q_1$ and
nondecreasing with respect to $q_2$.}
\\  ($\mathbf{G2}$)  \emph{consistency:
$g_j\left(x,q,q\right)=H_j(x,q)$  $\forall x\in [0,l_j]$, $\forall q \in \R $.}
\\($\mathbf{G3}$) \emph{differentiability: $g_j$ is  of class $\cC^1$.}
\\($\mathbf{G4}$) \emph{superlinear growth :
$ g_j(x,q_1,q_2)\ge \a((q_1^-)^2+(q_2^+)^2)^{\g/2}-C$ for some $\a>0$, $C\in\R$, $\g>1$  and $q^\pm_s$ denote the positive and negative part of $q_s$, $s=1,2$}
\\($\mathbf{G5}$)  \emph{convexity  : for all $x\in e_j$, $(q_1,q_2)\mapsto g_j\left(x,q_1,q_2\right)$ is convex.}\\
Numerical Hamiltonians fulfilling these requirements are provided by Lax-Friedrichs or Godunov type schemes,
see \cite{se}. As an example, suppose that the Hamiltonian $H$ is of the form
$H(x,p)= \Psi(x, |p|)$ where $\Psi$ is convex,  increasing and superlinear  with respect to its second
argument. Then the Engquist-Osher Godunov scheme reads as
\[
  g_j(x, q_1,q_2)=\Psi\big(x, ( \min(q_1,0)^2 +\max(q_2,0))^2 \big)
\]
and the  monotonicity, consistency and coercivity conditions  are satisfied.\\

Given $U,W:\cG_h\to \R$, we define the scalar product
\[(U,W)_2=\sum_{j\in J}\sum_{k=1}^{N^h_j-1}h_jU_{j,k} W_{j,k}+
          \sum_{i\in I} \left(\sum_{j\in \Inc^{+}_i}\frac{h_j}{2} U_{j,0} W_{j,0}+
          \sum_{j\in \Inc^{-}_i}\frac{h_j}{2}U_{j,N^h_j} W_{j,N^h_j}\right).\]
We introduce the compact and convex set
\[
\cK_h=\{ (M_{j,k})_{ j\in J,\,0\le k \le N^h_j}: \text{$M$ is continuous, } M_{j,k}\ge 0,\,  (M,1)_2=1 \}.
\]

The operator  $V[m](x_{j,k})$ is approximated by $(V_h[M])_{j,k}$ where $M$ is the piecewise constant function
taking the value $M_{j,k}$ in the interval  $|y-y_{j,k}|\le h_j/2$, $k=1,\dots,N^h_j-1$, $j\in J$
(at the vertices only the  half interval contained in $[0,l_j]$ is considered).
In particular, if $V$ is a local operator, i.e. $V[m](x)=F(m(x))$, then we set $(V_h[M])_{j,k}=F(M_{j,k})$.
We assume that\\\\
($\mathbf{V1}$) \emph{$V_h$ is continuous and maps   $\cK_h$  on a bounded set of grid functions.\\}
($\mathbf{V2}$)  \emph{$V_h$ is monotone, i.e.
 \[
    \left( V_h[M]-V_h[\bar M],M-\bar  M \right)_2 \le 0 \Rightarrow  M=  \bar M.
 \]}
($\mathbf{V3}$)
\emph{There exists $C$ independent of $h$ such that
 for all grid functions  $M\in \cK_h$
 \begin{align*}
 & \| V_h[M] \|_{\infty}:=\max_{j,k}|(V_h[M])_{j,k} | \le C\\
 & |(V_h[M])_{j,k} - (V_h[M])_{j,\ell} | \le C |y_{j,k}-y_{j,\ell}|\qquad k,\ell=0,\dots,N^h_j,\,j\in J.
 \end{align*}}
\section{A finite difference scheme for the stationary MFG system}\label{sec2}
In this section we introduce the approximation scheme for the system    \eqref{MFGs}.
For simplicity, we consider a network $\G$ without boundary; appropriate boundary condition can be inserted
in the scheme in a straightforward way.
 At the internal grid points we consider the finite difference  system
\begin{equation}\label{MFGDs}
\left\{
\begin{array}{ll}
    -\nu_j (D^2_h U)_{j,k} + g(x_{j,k}, [D_h U ]_{j,k} )+\Lambda = \left(V_h[M]\right)_{j,k},\quad& k=1,\dots, N^h_j-1,\, j\in J\\[8pt]
     \nu_j(D^2_h M)_{j,k}+\cB^h(U,M)_{j,k}=0,& k=1,\dots, N^h_j-1,\, j\in J\\[8pt]
     M\in \cK_h,\qquad (U,1)_2=0,
\end{array}
\right.
\end{equation}
\textbf{where $U$, $M$ are grid functions and $\Lambda\in\R$.}
The transport operator  $\cB^h $ is defined for $j\in J$ and  $k=1$ by
\[
\cB^h(U,M)_{j,k}=
 \begin{array}{l}
   \frac{1}{h_j}\left[M_{j,k}\frac{\pd g}{\pd q_1}(x_{j,k}, [D_h U]_{j,k}) +\right.\\[4pt]
  \left.    M_{j,k+1}\frac{\pd g}{\pd q_2}(x_{j,k+1}, [D_h U]_{j,k+1})-M_{j,k}\frac{\pd g}{\pd q_2}(x_{j,k}, [D_h U]_{j,k})\right];
 \end{array}
\]
for $k=2,\dots,N^h_j-2$ by
\[
    \cB^h(U,M)_{j,k}=
            \begin{array}{l}
                    \frac{1}{h_j}\left[M_{j,k}\frac{\pd g}{\pd q_1}(x_{j,k}, [D_h U]_{j,k})-M_{j,k-1}\frac{\pd g}{\pd q_1}(x_{j,k-1}, [D_h U]_{j,k-1})\right. \\[4pt]
                  +\left.M_{j,k+1}\frac{\pd g}{\pd q_2}(x_{j,k+1}, [D_h U]_{j,k+1})-M_{j,k}\frac{\pd g}{\pd q_2}(x_{j,k}, [D_h U]_{j,k})\right];
           \end{array}
\]
for  $k=N^h_j-1$ by
\[
\cB^h(U,M)_{j,k}=
 \begin{array}{l}
    \frac{1}{h_j}\left[M_{j,k}\frac{\pd g}{\pd q_1}(x_{j,k}, [D_h U]_{j,k})-M_{j,k-1}\frac{\pd g}{\pd q_1}(x_{j,k-1}, [D_h U]_{j,k-1})-\right.\\[4pt]
   \left. M_{j,k}\frac{\pd g}{\pd q_2}(x_{j,k}, [D_h U]_{j,k})\right].
 \end{array}
\]
We discuss now the transition conditions at the vertices, see \eqref{kir}. We  discretize the Kirchhoff
condition  for the function $u$ via a $1^{st}$ order approximation of the derivative and we impose
the continuity at the vertices
\begin{equation}\label{kirdiscru}
\left\{\begin{array}{ll}
   \cS^h(U, V_h[M]-\Lambda)_i=0,\qquad & i\in I, \\[4pt]
 \text{$U$ continuous at $v_i$},& i\in I,
\end{array}
\right.
\end{equation}
where for grid functions $U,V$, the operator $\cS^h:\cV\to \R$ is defined   by
\begin{equation}\label{transopu}
    \cS^h(U, V)_i= \sum_{j\in \Inc^{+}_i}  \big[\nu_j (D^+U)_{j,0}+\frac{h_j}{2} V_{j,0} \big]- \sum_{j\in \Inc^{-}_i}\big[ \nu_j (D^+ U )_{j,N^h_j-1}-\frac{h_j}{2}V_{j,N^h_j}\big] .
\end{equation}
To discretize the transition condition for $m$  we consider a $1^{st}$ order approximation of
the derivative (the continuity of $M$ at the vertices is included  in the definition of $\cK_h$)
\begin{equation}\label{kirdiscrm}
    \cT^h(M,U)_i=0\qquad i\in I,
\end{equation}
where for   grid functions $U,M$, the operator $\cT:\cV\to \R$ is defined  by
\begin{equation}\label{transopm}
\begin{split}
\cT^h(M,U)_i=& \sum_{j\in \Inc^{+}_i} \big[\nu_j (D^+M)_{j,0}+ M_{j,1}\frac{\pd g}{\pd q_2}(x_{j,1}, [D_h U]_{j,1})\big] \\
-& \sum_{j\in \Inc^{-}_i} \big[\nu_j (D^+ M )_{j,N^h_j-1}+M_{j,N^h_j-1}\frac{\pd g}{\pd q_1}(x_{j,N^h_j-1}, [D_h U]_{j,N^h_j-1})\big]=0.\\[4pt]
\end{split}
\end{equation}
\begin{Remark}
For the  discretization of the differential equations in \eqref{MFGs} inside the edge, we follow the
same approach in \cite{a,acd} and we refer to these papers for motivations and explanations. \textbf{We just recall
that the transport  operator $\cB^h$ comes from the discretization of the quantity
\[\ds \int_{e_j} m   H_p(x,\pd  u)  \pd  w\, dx\]
for a  test function $w$, which is connected with the weak formulation of the Fokker-Planck equation on the network.}\par
\textbf{ For the the approximation of the transition conditions in \eqref{kir},
 we use a standard $1^{st}$  order discretization  of the  normal derivative of $u$ and $m$ with the  sign depending
if the vertex corresponds to  either the initial point or the terminal one in the parametrization of the edge. The flux term in
the Kirchhoff condition for $m$ is approximated in a upwind fashion depending always on the orientation of the edge. Finally  the additional term
$\frac{h_j}{2}(\left(V_h[M]\right)-\Lambda\big)$ in \eqref{kirdiscru}, which vanishes for $h\to 0$, is  necessary  to obtain
the identity \eqref{fundid} which plays a key role in the uniqueness and convergence results. \\ }
\textbf{Note  that at  a vertex $v_i$, we have  respectively   $\#(\Inc_i)$    values $U_j$ and $\#(\Inc_i)$ values $M_j$, corresponding
to the restrictions of these functions to the incident edges $e_j$, $j\in \Inc_i$.
Since  \eqref{kirdiscru} and \eqref{kirdiscrm}  gives $\#(\Inc_i)$ linear conditions,
the value of $U$ and $M$  at $v_i$  is  univocally determined.}
\end{Remark}
Summarizing the approximation scheme for the stationary problem \eqref{MFGs} is given by
the \eqref{MFGDs}-\eqref{transopm}.
In the next subsections we study existence, uniqueness  and convergence of the scheme.
\subsection{Existence}
\textbf{We prove existence of a solution to \eqref{MFGDs}-\eqref{transopm}   by a fixed point argument.
We preliminarily need to prove      existence, uniqueness and  regularity   for the first equation in \eqref{MFGDs} with transition conditions \eqref{kirdiscru}, see Lemma \ref{th_ergodic}. This result is obtained, as in the continuous case, by approximating the limit ergodic problem \eqref{HJBs}  with  the sequence of problems \eqref{HJBerg}, which contains a zero order term $\rho U^\rho$,  and passing to the limit for $\rho\to 0$. For this we need to estimate, uniformly in $\rho$, the discrete gradient of $U^\rho$  (see Lemma \ref{erg2}). }
\begin{Lemma}
Let $V:\cG_h\to\R $ be a continuous grid function and assume that  $g$ satisfies \textbf{(G1)-(G3)}. For $\rho>0$,  there is a unique solution to the problem
\begin{equation}\label{HJBerg}
\left\{
\begin{array}{ll}
     -\nu_j (D^2_h U^\rho)_{j,k} + g(x_{j,k}, [D_h U^\rho ]_{j,k} )+\rho U^\rho_{j,k} =  V_{j,k},\qquad& k=1,\dots, N^h_j-1,\, j\in J\\[8pt]
     \cS^h(U^\rho, V- \rho U^\rho)_i=0,& i\in I;   \\[6pt]
     \text{$U^\rho$ continuous at $v_i$},& i\in I.
\end{array}
\right.
\end{equation}
\end{Lemma}
\begin{Proof}
To prove the existence we show that the map $\cF:\R^{N^h}\to\R^{N^h}$ defined by
\[
\cF(U)=\left\{
      \begin{array}{l}
       \frac{1}{\rho}\Big( \nu_j (D^2_h U )_{j,k} - g(x_{j,k}, [D_h U  ]_{j,k} )+V_{j,k}\Big),\quad \hbox{$j\in J$, $k=1,\dots,N^h_j-1$;} \\[8pt]
       \frac{1}{\rho h_{v_i}}\cS^h(U, V)_i, \quad \hbox{$i\in I$;} \\
      \end{array}
    \right.
\]
(where $h_{v_i}$ as in \eqref{hi}) admits a fixed point.\\
 Set $r=(\max_{j,k}|H(x_{j,k},0)|+\|V\|_\infty)/\rho$. By the regularity of $g$ the map $\cF$ is continuous from $B_r=\{U\in \R^{N^h}:\|U\|_\infty\le r\}$ to $ \R^{N^h}$.
Assume that $U\in \partial B_r$, hence $\max_{j\in J,\,k=0,\dots,N^h_j} |U_{j,k}| =r$. Consider first the case $U_{ j,  k}= r$ for some
 $ j\in J$, $ k\in \{1, \dots, N^h_j-1\}$.
Since $(D^2_hU)_{j,k}\le 0$, $D^+U_{j,k}\le 0$ and $D^+U_{j,k-1}\ge 0$, by the monotonicity and the consistency of $g$  we get
\[
\nu_j (D^2_h U)_{ j,  k} - g(x_{ j,  k}, [D_h U ]_{j,k} )\le -H(x_{ j,  k},0)
\]
and therefore
\[
\cF(U)_{ j,  k}\le \frac 1 \rho \left(-H(x_{ j,  k},0)+V_{ j,  k} \right)\le r.
\]
Hence $\cF(U)_{ j,  k}\le U_{ j,  k}$ and $\cF(U)_{ j,  k}\neq \mu U_{ j,  k}$ if $\mu>1$.\\
Now assume that there exists $i\in I$ such that  $U_i=r$ for some $i\in I$ ($U_i$ is the common value of $U_{j,k}$ at $v_i$)   then $(D^+U )_{j,0}\le 0$ if $v_i=\pi_j(0)$,
$(D^+U )_{j,N^h_j-1}\ge 0$ if $v_i=\pi_j(N^h_j)$
and therefore
\[
\cF(U)_i \le   \frac{2}{\rho h_{v_i}} \Big(   \sum_{j\in \Inc^{+}_i} \frac{h_j}{2}V_{j,0}+\sum_{j\in \Inc^{-}_i} \frac{h_j}{2}V_{j,N^h_j}\Big)\le r.
\]
Hence $\cF(U)_i\le U_i$ and $\cF(U)_i\neq \mu U_i$ if $\mu>1$.
Arguing in a similar way if either  $U_{ j,  k}=-r$ or $U_i=-r$,  we have that $\cF(U)\neq \mu U$ for all $\mu>1$ and $U\in \partial B_r$. Hence by the Leray-Schauder fixed point theorem there exists $U^\rho\in B_r$ such that $\cF(U^\rho)=U^\rho$ and therefore a solution of \eqref{HJBerg}. We also have the estimate
 \begin{equation}\label{erg1}
 \|U^\rho\|_\infty\le \frac{1}{\rho}(\max_{j,k}|H(x_{j,k},0)|+\|V\|_\infty).
 \end{equation}
We prove  uniqueness of the solution to \eqref{HJBerg}. Let   $U^1$, $U^2$ be  two solutions of \eqref{HJBerg} and assume by contradiction that
 $\max_{j,k}(U^1_{ j,k}-U^2_{j,k})=\d>0$. Consider  first the case  that there exists
 $\bar j\in J$, $\bar k\in \{1, \dots, N^h_j-1\}$ such that
 $U^1_{\bar j, \bar k}-U^2_{\bar j, \bar k}=\d $.
Subtracting the equations satisfied by $U^1$ and $U^2$, we get
\[-\nu_j (D^2_h (U^1-U^2))_{\bar j, \bar k} +g(x_{\bar j, \bar k}, [D_h U^1 ]_{j,k} )-g(x_{\bar j, \bar k}, [D_h U^2 ]_{j,k} )+\rho(U^1-U^2)_{\bar j, \bar k}=0.\]
Since $(\bar j, \bar k)$ is a maximum point for $U^1-U^2$, by the monotonicity of $g$, we get
\[\rho \d= \rho(U^1-U^2)_{\bar j, \bar k}\le 0\]
 and therefore a contradiction.
If there exists $i\in I$ such that $U^1_i-U^2_i=\d$, then subtracting  the transition conditions  satisfied by $U^1$ and $U^2$, we get
\begin{align*}
0=&\sum_{j\in \Inc^{+}_i}  ( \nu_j (D^+(U^1-U^2) )_{j,0}-\frac{h_j}{2}(\rho (U^1-U^2)_{j,0}) \\
 -&\sum_{j\in \Inc^{-}_i}( \nu_j (D^+(U^1-U^2) )_{j,N^h_j-1}
-\frac{h_j}{2}\rho (U^1-U^2)_{j,N^h_j})\le -\frac{\d\rho}{2} \sum_{j\in \Inc_i} h_j
\end{align*}
and therefore also in this case a contradiction. We conclude that $U^1\le U^2$ and we prove  in a similar way that $U^2\le U^1$.
\end{Proof}
In the next lemma, we get an a priori bound for  the  gradient of the solution to the discrete Hamilton-Jacobi-Bellman equation
by assuming that the function is bounded. It is important for the analysis of the convergence of the scheme that all the bounds are uniform in $h$.
\begin{Lemma}\label{th_apriori}
Let $\tilde V:\cG_h\to\R $ be a  continuous grid function and assume that  $g$ satisfies \textbf{(G1)-(G4)}.
Let $U^h$ be a solution of the problem
\begin{equation}\label{apr1}
\left\{
\begin{array}{ll}
    -\nu_j (D^2_h U)_{j,k} + g(x_{j,k}, [D_h U ]_{j,k} ) =  \tilde V_{j,k},\qquad  &k=1,\dots, N^h_j-1,\, j\in J\\[8pt]
\cS^h(U,\tilde V)_i=0  &i\in I;   \\[6pt]
     \text{$U$ continuous at $v_i$},& i\in I,
\end{array}
\right.
\end{equation}
and assume that
\begin{equation}\label{apr2}
\|U^h\|_\infty\le C_0
\end{equation}
with $C_0$ independent of $h$.
Then
\[
 \|D_h  U^h\|_\infty:=\max_{j\in J}\max_{k=0,\dots, N^h_j-1}|(D^+ U^h)_{j,k}|\le C
\]
where $C$ depends on $C_0$, $\|\tilde V\|_\infty$, but not on   $h$.
\end{Lemma}
\begin{Proof}
We first prove that $D^+ U^h$ is bounded at the vertices. Assume by contradiction that for some $i\in I$
\[\max\left\{\max_{j\in \Inc^{+}_i}   |(D^+U^h)_{j,0}|, \max_{j\in \Inc^{-}_i}   |(D^+U^h)_{j,N^h_j-1}|\right\}\to+ \infty\quad \text{for $|h|\to 0$}.\]
Because of the transition condition in  \eqref{apr1}, it is not restrictive to assume that, up to a subsequence,
\[\max\left\{\max_{j\in \Inc^{+}_i}    \{(D^+U^h)_{j,0}\} , \max_{j\in \Inc^{-}_i} \{ - (D^+U^h)_{j,N^h_j-1}\}\right\}\to +\infty\quad \text{for $|h|\to 0$}.\]
Hence we assume that there exists $j\in \Inc^{+}_i$ such that $D^+U^h_{j,0}\to +\infty$ for $h_j\to 0$ (we proceed in a similar
way if there exists $j\in \Inc^{-}_i$ such that $-D^+U^h_{j,N^h_j-1}\to +\infty$).\\
Let $h_0$ be such that for $h_j<h_0$
\[
     D^+U^h_{j,0}\ge \frac{1}{\a}(C+\|  \tilde V\|_\infty)+ \frac{4C_0}{l_j}
\]
where $C$ as in \textbf{$(G4)$}, $C_0$ as in \eqref{apr2}.
Since $(D^2_h U^h)_{j,1}=(D^+U^h_{j,1}-D^+U^h_{j,0})/h$ we have
\begin{equation}\label{apr4a}
\begin{split}
   \frac {\nu_j}{h_j}D^+U^h_{j,1}=\frac {\nu_j}{h_j}D^+U^h_{j,0}+g(x_{j,1}, [D_h U^h ]_{j,1} ) -  \tilde V_{j,k}\ge\\
   \frac {\nu_j}{h_j}D^+U^h_{j,0}+ \a|D^+U^h_{j,0}|^\g-C-\|\tilde V\|_\infty\ge \frac{\nu_j}{h_j}D^+U^h_{j,0}
\end{split}
\end{equation}
and therefore $D^+U^h_{j,1}\ge D^+U^h_{j,0}$.
Iterating the previous inequality, we get
\begin{equation}\label{apr5}
   D^+U^h_{j,k+1}\ge D^+U^h_{j,k}  \quad\text{for $k=0,\dots,N^h_j-1$}.
\end{equation}
For $L\le N^h_j-1$, we have
\[\left\{\begin{array}{l}
    U^h_{j,1}=U^h_{j,0}+h_jD^+U^h_{j,0},\\
    U^h_{j,2}=U^h_{j,1}+h_jD^+U^h_{j,1}=U^h_{j,0}+h_j(D^+U^h_{j,0}+D^+U^h_{j,1}) ,\\
    \vdots\\
    U^h_{j,L}=U^h_{j,0}+h_j\sum_{k=0}^{L-1} D^+U^h_{j,k}.
\end{array}\right.\]
If $Lh_j > l_j/2$, by \eqref{apr5} we get
\[
U^h_{j,L}\ge U^h_{j,0}+Lh_j D^+ U^h_{j,0}\ge  U^h_{j,0}+ Lh_j \frac{4C_0}{l_j}> C_0
\]
and therefore a contradiction to \eqref{apr2}. \par
We show that $D^+ U^h$ is bounded also inside $\G$. Assume by contradiction that
there exists $j\in J$, $k_h\in\{1,\dots,N^h_j -2\}$ such that, up to a subsequence,
\begin{equation}\label{apr6}
    |D^+U^h_{j,k_h}|\to +\infty\qquad \text{for $h\to 0$.}
\end{equation}
By compactness, $x_{j,k_h}\to x_0\in e_j$ for $h\to 0$. We set  $y_0=\pi_j^{-1}(x_0)\in [0,l_j]$
and we  first consider the case $y_0\in (0,l_j)$.
If  $D^+U^h_{j,k_h}\to +\infty$ for $h_j\to 0$, let $h_0$ be such that
for $h_j\le h_0$
  \[
     D^+U^h_{j,k_h}\ge \frac{1}{\a}(C+\|\tilde V\|_\infty)+ \frac{4C_0}{l_j-y_0}.
\]
Arguing as in \eqref{apr4a}, we have
\[ D^+U^h_{j,k_h+l}\ge D^+U^h_{j,k_h+l-1}  \quad\text{for $l=1,\dots,N^h_j-k_h$}.\]
For $Lh_j> (l_j-y_0)/2$ we get
\[D^+U^h_{j,k_h+L}= U^h_{j,k_h}+h_j\sum_{l=0}^{L-1} D^+U^h_{j,k_h+l}\ge -C_0+Lh \frac{4C_0}{l_j-y_0}> C_0  \]
and therefore a contradiction to  \eqref{apr2}.\\
We now  consider the case $y_0\in (0,l_j)$ and $D^+U^h_{j,k_h} \to-\infty$ for $h\to 0$. Let $h_0$ be such that
for $h_j\le h_0$
\[
     D^+U^h_{j,k_h}\le- \frac{1}{\a}(C+\|\tilde V\|_\infty)-\frac{4C_0}{y_0}
\]
We have
\[\begin{split}
   \frac {\nu_j}{h_j}D^+U^h_{j,k_h}=\frac {\nu_j}{h_j}D^+U^h_{j,k_h-1}+g(x_{j,1}, [D_h U^h ]_{j,k_h} ) -  \tilde V^h_{j,k_h}\ge\\
   \frac {\nu_j}{h_j}D^+U^h_{j,k_h-1}+ \a|D^+U^h_{j,k_h}|^\g-C-\|\tilde V\|_\infty\ge \frac{\nu_j}{h_j}D^+U^h_{j,k_h-1}
\end{split}\]
and iterating
\begin{equation}\label{apr5a}
   D^+U^h_{j,k_h-l}\ge D^+U^h_{j,k_h-l-1}  \quad\text{for $l=0,\dots,k_h-1$}.
\end{equation}
For $L\le k_h-1$, we have
\[
\left\{\begin{array}{l}
    U^h_{j,k_h-1}=U^h_{j,k_h}-hD^+U^h_{j,k_h-1},\\
    \vdots\\
    U^h_{j,k_h-L}=U^h_{j,k_h}-h\sum_{l=1}^{L} D^+U^h_{j,k_h-l}.
\end{array}\right.
\]
Hence if  $Lh_j > y_0/2$, by \eqref{apr5a} we get
\[
U^h_{j,k_h-L}\ge U^h_{j,k_h}-h_j\sum_{l=1}^{L} D^+U^h_{j,k_h-l}\ge U^h_{j,k_h}-Lh_j D^+U^h_{j,k_h}\ge  -C_0+ Lh_j \frac{4C_0}{y_0}> C_0
\]
and therefore a contradiction to \eqref{apr2}. \par
In case $y_0=\pi_j^{-1}(x_0)$ is equal  either to $0$ or  to $l_j$ and $|D^+U^h_{j,k_h}|\to +\infty$, it is easy to adapt the previous arguments to obtain again
a contradiction to \eqref{apr2}.
\end{Proof}
\begin{Lemma}
Let $U^\rho$ be the solution of \eqref{HJBerg}, then
\begin{equation}\label{erg2}
 \|D^+U^\rho\|_\infty \le C_2
\end{equation}
for a constant $C_2$ independent  of $\rho$ and $h$.
\end{Lemma}
\begin{Proof}
Fix an arbitrary node $x_{\bar j,\bar k}\in \G$ and set $W^\rho=U^\rho-U^\rho_{\bar j,\bar k}$.
Adapting to the case of the networks the argument  in \cite[Prop.2]{acd}, it is possible to show  that $W^\rho$ is bounded, uniformly in $\rho$.
Since $W^\rho$ is a solution of \eqref{apr1} with $\tilde V=V-\rho U^\rho$ and by \eqref{erg1}
$\tilde V$ is bounded, uniformly in $\rho$ and $h$, we can apply Lemma \ref{th_apriori} to get a bound
on $\|DW^\rho\|_\infty$ and therefore on $\|DU^\rho\|$ uniform  in $\rho$ and $h$.
\end{Proof}
\begin{Lemma}\label{th_ergodic}
Let $V:\cG_h\to\R $ be a continuous grid function and assume that  $g$ satisfies \textbf{(G1)-(G4)}. Then there exists a unique
couple   $(U,\Lambda)$, where $U:\cG_h\to\R$ and $\Lambda\in\R$, solution of the problem
\begin{equation}\label{HJBs}
\left\{
\begin{array}{ll}
    -\nu_j (D^2_h U)_{j,k} + g(x_{j,k}, [D_h U ]_{j,k} )+\Lambda =  V_{j,k},\qquad  &k=1,\dots, N^h_j-1,\, j\in J\\[8pt]
\cS^h(U,V-\Lambda)_i=0,&i\in I;   \\[6pt]
     \text{$U$ continuous at $v_i$},&i\in I\\[4pt]
     (U,1)_2=0.
\end{array}
\right.
\end{equation}
Moreover
\begin{equation}\label{erg}
  |\Lambda|\le C_1, \quad  \|D_h U\|_\infty\le C_2
\end{equation}
for some constants  $C_1$, $C_2$   independent of  $h$.
\end{Lemma}
\begin{Proof}
We prove existence by passing to the limit in the ergodic approximation \eqref{HJBerg}. By \eqref{erg1}
\[
   \|\rho U^\rho\|_\infty\le C_1
\]
for any $\rho>0$ where $C_1$ is independent of $\rho$.
By \eqref{erg1} and  \eqref{erg2}, up to a subsequence,   $U^\rho-(U^\rho,1)_2$ converges to a function $U:\cG_h\to\R$ such that $(U,1)_2=0$ and  $\rho U^{\rho}_{j,k}$ converges to  $\Lambda\in\R$
(independent of $(j,k)$). Moreover  the couple $(U,\Lambda)$ satisfies \eqref{HJBs} and  the bounds in \eqref{erg}.\\
The uniqueness of the  couple $(U,\Lambda)$ can be proved by an argument similar to the one for the uniqueness of \eqref{HJBerg}.
\end{Proof}
\begin{Remark}
Note  that the dependence of the  bounds in  \eqref{erg} on the function $V_{j,k} $ is only
by means of $\|V\|_\infty$. This is crucial for the proof of the next theorem.
\end{Remark}

\begin{Theorem}\label{th_existence}
If  $g$ satisfies \textbf{(G1)-(G4)}, $V$ satisfies \textbf{(V1)}, then the problem  \eqref{MFGDs}-\eqref{transopm} has at least a solution  $(U,M,\Lambda)$. Moreover
\begin{equation}\label{stimasol}
  |\Lambda|\le C_1,\qquad  \| U\|_\infty+ \|D_h U\|_\infty\le C_2
\end{equation}
for some constants  $C_1$, $C_2$ independent of $h$.
\end{Theorem}
\begin{Proof}
We define a map $\Phi$ which associates to $M\in \cK_h$  the solution $(U,\Lambda)$
of the  problem \eqref{HJBs} with $V_{j,k}=(V_h[M])_{j,k}$. By   Lemma \ref{th_ergodic}, the map $\Phi$ is well defined. \\
We show that $\Phi$ is continuous.
Let $M^s\in \cK_h$ be such that $M^s\to M\in \cK_h$ as $s\to \infty$, hence by \textbf{(V1)}, $V_h[M^s]\to V_h[M]$ as $s\to \infty$. Let $(U^s, \Lambda^s)$ be the sequence of solutions of \eqref{HJBs} with $V=V_h[M^s]$. By \eqref{erg} the sequences $\Lambda^s$ and $\|U^s\|_\infty$ are bounded and therefore, up to a subsequence, converge to $\Lambda \in\R$ and, respectively, to a grid function $U$. It is immediate that $(U,\Lambda)$ is a solution of \eqref{HJBs} with $V=V_h[M]$.
By the uniqueness of the solution to \eqref{HJBs}, it follows that all the sequence $(\Lambda^s, U^s)$ converges to $(\Lambda,U)$ and therefore the continuity of the map $\Phi$ and the estimate \eqref{stimasol}.\\
We define a map $\Psi$ which associates to $M\in \cK_h$ the solution $\bar M$ of linear problem
\begin{equation*}
\left\{
\begin{array}{l}
  \mu  \bar M_{j,k}-\nu_j(D^2_h\bar M)_{j,k}-\cB^h(U,\bar M)_{j,k}=\mu M_{j,k}\qquad   j\in J,\,k=1,\dots,N^h_j-1\\[6pt]
  \mu \bar M_{i}- \sum_{j\in \Inc^{+}_i}\frac{2}{h_j} \big[\nu_j (D^+\bar M)_{j,0}+ \bar M_{j,1}\frac{\pd g}{\pd q_2}(x_{j,1}, [D_h U]_{j,1})\big]+ \\
\sum_{j\in \Inc^{-}_i}\frac{2}{h_j} [\nu_j (D^+\bar  M )_{j,N^h_j-1}-\bar M_{j,N^h_j-1}\frac{\pd g}{\pd q_1}(x_{j,N^h_j-1}, [D_h U]_{j,N^h_j-1})]= \mu M_{i} \qquad i\in I
\end{array}
\right.
\end{equation*}
where $\mu>0$ and  $(U,\Lambda)=\Phi(M)$.
We rewrite the previous problem as
\begin{equation}\label{ex2}
    \mu \bar M+\cA\bar M=\mu M
\end{equation}
where $\cA$ is $N^h\times N^h$ matrix.
By the monotonicity and the regularity of $g$,  for $\mu$ sufficiently large the matrix $\mu I +\cA$ is a non singular $M$-matrix
 and is therefore invertible.
It follows that for any $M\in \cK_h$, \eqref{ex2} admits a solution $\bar M$  \textbf{and by $M$-matrix property  $\bar M\ge 0$ since $M\ge 0$.}
We prove that $(\bar M,1)_2=1$. First observe that if  $W,Z:\cG_h\to\R$, then
\begin{equation}\label{ex2b}
\begin{split}
  \sum_{j\in J}  \sum_{k=1}^{N^h_j-1}&\nu_j (D^2_h W)_{j,k}Z_{j,k}=- \sum_{j\in J}\sum_{k=1}^{N^h_j-2}\nu_j (D^+W)_{j,k}(D^+Z)_{j,k}
    \\&-\sum_{i\in I}\sum_{j\in \Inc^{+}_i}\frac{\nu_j}{h_j} Z_{j,1} (D^+  W)_{j,0}
    +\sum_{i\in I}\sum_{j\in \Inc^{-}_i}\frac{\nu_j}{h_j} Z_{j,N^h_{j-1}} (D^+  W)_{j,N^h_j-1}
\end{split}
\end{equation}
and
\begin{equation}\label{ex2c}
\begin{split}
\sum_{j\in J}\sum_{k=1}^{N^h_j-1}& \cB^h(U,W)_{j,k} Z_{j,k}  =-\sum_{j\in J}\sum_{k=1}^{N^h_j-1}
W_{j,k}\,[D_h  Z]_{j,k}\cdot \nabla_q g(x_{j,k},[D_hU]_{j,k})\\
&-\sum_{i\in I}\Big(\sum_{j\in \Inc^{+}_i}\frac{1}{h_j}\big[W_{j,1}Z_{j,0}\frac{\pd g}{\pd q_2}(x_{j,1},[D_hU]_{j,1})\big]\\
&-\sum_{j\in \Inc^{-}_i} \frac{1}{h_j}\big[W_{j,N^h_j-1}Z_{j,N^h_j}\frac{\pd g}{\pd q_1}(x_{j,N^h_j-1},[D_hU]_{j,N^h_j-1})\big]\Big)
 \end{split}
\end{equation}
If $W=\bar M$, $Z\equiv 1$, by \eqref{ex2b}-\eqref{ex2c} we get
\begin{align*}
& \sum_{j\in J}\sum_{k=1}^{N^h_j-1} (\cA \bar M)_{j,k}  =  \sum_{i\in I}\sum_{j\in \Inc^{+}_i}\frac{1}{h_j}\Big[\nu_j(D^+  \bar M)_{j,0}+  \bar M_{j,1} \frac{\pd g}{\pd q_2}(x_{j,1},[D_hU]_{j,1})\Big]\\
&-\sum_{i\in I}\sum_{j\in \Inc^{-}_i} \frac{1}{h_j}\Big[\nu_j(D^+  \bar M)_{j,N^h_j-1}-\bar M_{j,N^h_j-1} \frac{\pd g}{\pd q_1}(x_{j,N^h_j-1},[D_hU]_{j,N^h_j-1})\Big]
\end{align*}
Hence by the definition of $\cA$ at the vertices we have
\[
\sum_{j\in J}\sum_{k=1}^{N^h_j-1} (\cA \bar M)_{j,k}=\frac{\mu} 2 \sum_{i\in I}\left(\sum_{j\in Inc^{+}_i}   (\bar M-M)_{j,0}+  \sum_{j\in \Inc^{-}_i}  (\bar M-M)_{j,N^h_j}\right).
\]
Therefore
\[  \sum_{j\in J}\sum_{k=1}^{N^h_j-1}h_j [\mu \bar M_{j,k} +(\cA \bar M)_{j,k} ]=\mu\sum_{j\in J}\sum_{k=1}^{N^h_j-1}h_j     M_{j,k}  \]
which implies  $(\bar M,1)_2=(M,1)_2=1$ and therefore $\bar M\in\cK_h$.\\
Hence    $\Psi$ maps $\cK_h$ into $\cK_h$. From the boundedness and continuity of $\Phi$
and the regularity of $g$,   $\Psi$ is continuous. By the Brouwer's fixed point theorem it follows that $\Psi$ admits a fixed point $M$  which is a solution of   \eqref{MFGDs}-\eqref{transopm}.
\end{Proof}
\subsection{Uniqueness}
We first prove a fundamental identity  which plays a crucial role in  uniqueness and convergence of the scheme (compare with
\cite[(3.20)]{acc}.
\begin{Lemma}
Let $A, B:\cG_h\to\R$ be two   grid functions, $(U,M,\Lambda)$   a solution of \eqref{MFGDs}-\eqref{transopm} and $(\bar U,\bar M,\bar \Lambda)$ a   solution  of
\begin{equation}\label{MFGDpert}
\left\{
\begin{array}{ll}
    -\nu_j (D^2_h \bar U)_{j,k} + g(x_{j,k}, [D_h\bar U ]_{j,k} )+\bar\Lambda = \left(V_h[\bar M]\right)_{j,k}+A_{j,k},\quad& k=1,\dots, N^h_j-1,\, j\in J\\[4pt]
     \nu_j(D^2_h\bar M)_{j,k}+\cB^h(\bar U,\bar M)_{j,k}=B_{j,k},& k=1,\dots, N^h_j-1,\, j\in J\\[4pt]
     \cS^h(\bar U,V_h[\bar M]-\Lambda)_i=\sum_{j\in \Inc^+_i}\frac{h_j}{2}A_{ j,0 }+\sum_{j\in \Inc^-_i}\frac{h_j}{2}A_{ j,N^h_j }, &i\in I\\[4pt]
     \cT^h(\bar M,\bar U)_i=\sum_{j\in \Inc^+_i}\frac{h_j}{2}B_{ j,0 }+\sum_{j\in \Inc^-_i}\frac{h_j}{2}B_{ j,N^h_j },&i\in I\\[4pt]
\bar U \textrm{continuous at } v_i,&i\in I\\[4pt]
    \bar M\in \cK_h,\qquad (\bar U,1)_2=0.
\end{array}
\right.
\end{equation}
Then
\begin{equation}\label{fundid}
   \cR^h(M,U,\bar U)+\cR^h(\bar M,\bar U,  U)+(V_h[M]-V_h[\bar M], M-\bar M)_2=(A,M-\bar M)_2+(B,U-\bar U)_2
\end{equation}
where
\begin{align*}
    \cR^h(M,U,\bar U)=
\sum_{j\in J}\sum_{k=1}^{N^h_j-1}
 h_j &M_{j,k}\Big[
g(x_{j,k}, [D_h\bar U]_{j,k})-g(x_{j,k}, [D_h  U]_{j,k})\\
&-[D_h (\bar U- U)]_{k,j}\cdot \nabla_q g(x_{j,k}, [D_h U]_{j,k})\Big]
\end{align*}
\end{Lemma}
\begin{Proof}
Let $(U,M,\Lambda)$ and $(\bar U,\bar M,\bar \Lambda)$ be as in the statement. 
Subtracting the equations for $U$ and $\bar U$,  multiplying the resulting equation by $h_j(M-\bar M)_{j,k}$
and summing over $j\in J$, $k=1,\dots,N^h_j-1$, we get
\begin{equation}\label{un1}
\begin{split}
    \sum_{j\in J}\sum_{k=1}^{N^h_j-1}& h_j (M-\bar M)_{j,k}\Big[-\nu_j D^2_h(U-\bar U)_{j,k}+ g(x_{j,k}, [D_h U]_{j,k})-g(x_{j,k}, [D_h \bar U]_{j,k})+\\
     &(\Lambda-\bar\Lambda)-(V_h[M]-V_h[\bar M])_{j,k}\Big]=\sum_{j\in J}\sum_{k=1}^{N^h_j-1} h_j(M-\bar M)_{j,k}A_{j,k}.
\end{split}
\end{equation}
Subtracting the equations  for   $M$ and $\bar M$, multiplying the resulting equation by  $h_j(U-\bar U)_{j,k}$ and summing over $j\in J$, $k=1,\dots,N^h_j-1$, we get
\begin{equation} \label{un2}
\begin{split}
    \sum_{j\in J}\sum_{k=1}^{N^h_j-1}& h_j (U-\bar U)_{j,k}\left[\nu_j D^2_h (M-\bar M)_{j,k}    +\cB^h(U,M)_{j,k}-\cB^h(\bar U, \bar M)_{j,k}\right]   =\\
  & \sum_{j\in J}\sum_{k=1}^{N^h_j-1} h_j(U-\bar U)_{j,k}B_{j,k}.
\end{split}
\end{equation}
We have the identity
\begin{equation}\label{un3}
\begin{split}
   -\sum_{j\in J}\sum_{k=1}^{N^h_j-1}  \nu_j (U-\bar U)_{j,k} D^2_h (M-\bar M)_{j,k}=\sum_{j\in J}\sum_{k=1}^{N^h_j-1} -\nu_j (M-\bar M)_{j,k} D^2_h  (U-\bar U)_{j,k}\\
  + \sum_{j\in \Inc^{+}_i}\frac{\nu_j}{h_j}\Big[ (M-\bar M)_{j,0}(D^+ (U-\bar U))_{j,0}- (U-\bar U)_{j,0}(D^+(M-\bar M))_{j,0}\Big]\\
  + \sum_{j\in \Inc^{-}_i} \frac{\nu_j}{h_j}\Big[(M-\bar M)_{j,N^h_j}(D^+ (U-\bar U))_{j,N^h_{j-1}}- (U-\bar U)_{j,N^h_j}(D^+(M-\bar M))_{j,N^h_{j-1}}\Big]
 \end{split}
\end{equation}
and, respectively,
\begin{equation}\label{un4}
\begin{split}
&\sum_{j\in J}\sum_{k=1}^{N^h_j-1} \cB^h(U,M)_{j,k} (U-\bar U)_{j,k}  =-\sum_{j\in J}\Big[\sum_{k=2}^{N^h_j-2}
M_{j,k}\,[D_h  (U-\bar U)]_{j,k}\cdot \nabla_q g(x_{j,k},[D_hU]_{j,k})\\
&-   M_{j,1}      \frac{\pd g}{\pd q_2}(x_{j,1}, [D_h U]_{j,1}) (\bar U-U)_{j,0} +M_{j,N^h_j-1} \frac{\pd g}{\pd q_1}(x_{j,N_j^h-1}, [D_h U]_{j,N_j^h-1}) (\bar U-U)_{j,N_j^h}\Big].
\end{split}
\end{equation}
In a similar way a  corresponding equation  for $\sum_{j\in J}\sum_{k=1}^{N^h_j-1} \cB^h(\bar U,\bar M)_{j,k} (U-\bar U)_{j,k}$ is also obtained. \\
We now discuss the boundary terms in \eqref{un3} and \eqref{un4}.
By the transition conditions for $U$ and $\bar U$ and the continuity  of $M$ at the vertices we have
\begin{equation}\label{un5}
\begin{split}
& -\sum_{j\in \Inc^{+}_i}\frac{\nu_j}{h_j} (M-\bar M)_{j,0}(D^+ (U-\bar U))_{j,0}+\sum_{j\in \Inc^{-}_i} \frac{\nu_j}{h_j} (M-\bar M)_{j,N^h_j}(D^+ (U-\bar U))_{j,N^h_{j-1}}\\
& =\sum_{j\in \Inc^{+}_i}  \frac{1}{2}(M-\bar M)_{j,0}[(V_h[M]-V_h[\bar M])_{j,0}-(\Lambda-\bar \Lambda)-A_{j,0} ]\\
 &   +   \sum_{j\in \Inc^{-}_i}  \frac{1}{2}(M-\bar M)_{j,N^h_j}[(V_h[M]-V_h[\bar M])_{j,N^h_j}-(\Lambda-\bar \Lambda)-A_{j,N_j^h}].
\end{split}
\end{equation}
By  transition conditions for $M$ and $\bar M$  we get
\begin{align}
\begin{split}\label{un6}
  &\sum_{j\in  \Inc^+_i} \frac{\nu_j}{h_j} (U-\bar U)_{j,0}   (D^+(M-\bar M))_{j,0}-\sum_{j\in \Inc^{-}_i} \frac{\nu_j}{h_j}  (U-\bar U)_{j,N^h_j}(D^+ (M-\bar M)  )_{j,N^h_j-1}\\
  &=-\sum_{j\in \Inc^+_i} \frac{(U-\bar U)_{j,0}}{h_j}\Big[M_{j,1}\frac{\pd g}{\pd q_2}(x_{j,1}, [D_h U]_{j,1})-\bar M_{j,1}\frac{\pd g}{\pd q_2}(x_{j,1}, [D_h \bar U]_{j,1})\Big]\\
&+\sum_{j\in \Inc^{-}_i}\frac{(U-\bar U)_{j,N^h_j}}{h_j}\Big[M_{j,N^h_j-1}\frac{\pd g}{\pd q_1}(x_{j,N^h_j-1}, [D_h U]_{j,N^h_j-1})  \\
&-\bar M_{j,N^h_j-1} \frac{\pd g}{\pd q_1}(x_{j,N^h_j-1}, [D_h \bar U]_{j,N^h_j-1})\Big]+ \sum_{j\in \Inc^+_i}\frac{1}{2}B_{ j,0 }(U-\bar U)_{j,0}\\
&+\sum_{j\in \Inc^-_i}\frac{1}{2}B_{ j,N^h_j }(U-\bar U)_{j,N^h_j}.
 \end{split}
\end{align}
Replacing \eqref{un3}--\eqref{un6} in \eqref{un2} and adding the resulting equation   to \eqref{un1} we finally get
(recall  that $(M,1)_2=(\bar M,1)_2=1$)
\begin{align*}
 \sum_{j\in J}\sum_{k=1}^{N^h_j-1}
 h_j\Big[&M_{j,k}
\left(
\begin{array}{c}
g(x_{j,k}, [D_h\bar U]_{j,k})-g(x_{j,k}, [D_h  U]_{j,k})\\
-[D_h (\bar U- U)]_{k,j}\cdot \nabla_q g(x_{j,k}, [D_hU]_{j,k})
\end{array}
\right)
\\ & +     \bar M_{j,k}\left(
\begin{array}{c} g(x_{j,k}, [D_h U]_{j,k})-g(x_{j,k}, [D_h  \bar U]_{j,k})\\
-[D_h ( U- \bar U)]_{k,j}\cdot \nabla_q g(x_{j,k}, [D_h \bar U]_{j,k}
\end{array}\right)
\Big]
\\[4pt]&+(V_h[M]-V_h[\bar M], M-\bar M)_2=(A,M-\bar M)_2+(B,U-\bar U)_2,
\end{align*}
which amounts to \eqref{fundid}.
\end{Proof}
\begin{Theorem}\label{th_uniq}
If $g$ satisfies  \textbf{(G1)-(G5)} and the operator $V_h$ is strictly monotone, i.e.
 \begin{equation}\label{strictlymon}
 (V_h[M]-V_h[\bar M],M-\bar M)_2 \le 0\,\Rightarrow M=\bar M
 \end{equation}
 then the problem \eqref{MFGDs}-\eqref{transopm} has at most one solution.
\end{Theorem}
\begin{Proof}
Let $(U,M,\Lambda)$ and $(\bar U,\bar M,\bar \Lambda)$ be  two solutions of  \eqref{MFGDs}-\eqref{transopm}. By \eqref{fundid}
with $A\equiv B\equiv 0$ we get
\[\cR^h(M,U,\bar U)+\cR^h(\bar M,\bar U,  U)+(V_h[M]-V_h[\bar M], M-\bar M)_2=0.\]
By the convexity of $g$ and the monotonicity of $V$, we see that all the terms in the left hand side of the previous equality are positive and therefore must vanish. The strong monotonicity of $V$ implies that $M=\bar M$. Hence   $(U,\Lambda)$, $(\bar U,\bar \Lambda)$ solve \eqref{HJBs}  with $V_{j,k}=V_h[M]_{j,k}=V_h[\bar M]_{j,k}$ and  by Lemma  \ref{th_ergodic} we get $U=\bar U$ and $\Lambda=\bar\Lambda$.
\end{Proof}
\subsection{Convergence}
In this section we analyze  the convergence of the scheme \eqref{MFGDs}-\eqref{transopm}  in the reference case
\begin{equation}\label{hamconv}
    H(x,p)=|p|^\b+f(x)
\end{equation}
where $\b\ge 2$ and $f:\G\to\R$ is a continuous function.  By \cite{cm}, we know that in this case
there exists a unique solution $(u,m,\l)$ to \eqref{MFGs} with $u\in C^{2,\a}(\G)$, $m\in C^2(\G)$, $m>0$,
and $\l\in\R$.

We consider a numerical Hamiltonian of the form
\begin{equation}\label{hamnumconv}
    g(x,p)=G(p_1^-,p_2^+)+f(x)
\end{equation}
  where  $G(p_1,p_2)=(p_1^2+p_2^2)^{\b/2}$ and $p^\pm_s$ denote the positive and negative part of $p_s$, $s=1,2$. We observe that $g$ satisfies  assumptions \textbf{(G1)-(G5)}.
We need   an   additional  assumption:\\[4pt]
($\mathbf{V4}$)  \emph{For any $m\in \cK:=\{\mu\in C^{0,\a}(\G):\, \int_\G \mu dx=1\}$, $M\in\cK^h$, denoted by $\cI^h(M)$ the continuous piecewise linear
reconstruction of $M\in \cK^h$ on $\G$, then
\[
\|V[m]-V_h[M]\|_\infty \le \o(\|m-\cI^h(M)\|_\infty)
\]
where $\o$ is a continuous, increasing function such that $\lim_{t\to 0^+} \o(t)=0$.}\\[4pt]
In the following we denote by  $o(1)$   a generic grid function whose maximum norm tends to $0$ as $|h|\to 0$.
Given a solution   $(u,m,\l)$ of \eqref{MFGs}, we define a grid function $u^h$ by
\begin{align*}
   & u^h_{j,k}:=\frac{1}{h_j}\int_{|y-y_{j,k}|\le h_j/2}u_j(y)dy,\quad \text{if $j\in J, k=1,\dots,N_j^h-1$},\\
   & u^h_i:=\sum_{j\in \Inc_i^+}\frac{2}{h_j}\int_{0\le y-y_{j,0}\le h_j/2}u_j(y)dy   +     \sum_{j\in \Inc_i^-}\frac{2}{h_j}\int_{0\le y_{j,N^h_j}-y\le h_j/2}u_j(y)dy,             \quad \text{if $i\in I$.}
\end{align*}
(note that $(u^h,1)_2=0$). We define in a similar way the grid function $m^h\in \cK_h$ and we also set $\l^h:=\l$.
Observe that by  ($\mathbf{V4}$)
\begin{equation}\label{conv0}
 \lim_{h\to 0} \|V[m]-V_h[m^h]\|_\infty =0.
\end{equation}
Hence by \eqref{conv0} and  the consistency assumption  ($\mathbf{G2}$),   $(u^h,m^h,\l^h)$ is a solution of
\begin{equation}\label{conv1a}
\left\{
\begin{array}{ll}
    -\nu_j (D^2_h u^h)_{j,k} + g(x_{j,k}, [D_h u^h ]_{j,k} )+\l^h = \left(V_h[m^h]\right)_{j,k}+A^h_{j,k},\quad& k=1,\dots, N^h_j-1,\, j\in J\\[8pt]
     \nu_j(D^2_h m^h)_{j,k}+\cB^h(u^h,m^h)_{j,k}= B^h_{j,k},& k=1,\dots, N^h_j-1,\, j\in J\\[8pt]
     m^h\in \cK_h,\qquad (u^h,1)_2=0
\end{array}
\right.
\end{equation}
with the transition conditions
\begin{equation}\label{conv1b}
\left\{\begin{array}{ll}
     \cS^h(u^h,V_h[m^h]-\l^h)_i=\sum_{j\in \Inc^+_i}\frac{h_j}{2}A^h_{ j,0 }+\sum_{j\in \Inc^-_i}\frac{h_j}{2}A^h_{ j,N^h_j }, &i\in I\\[4pt]
     \cT^h(m^h,u^h)_i=\sum_{j\in \Inc^+_i}\frac{h_j}{2}B^h_{ j,0 }+\sum_{j\in \Inc^-_i}\frac{h_j}{2}B^h_{ j,N^h_j },&i\in I\\[4pt]
       \text{$u^h$ continuous at $v_i$},&i\in I
\end{array}
\right.
\end{equation}
where $A^h$, $B^h$ are two grid functions such that
\begin{equation}\label{conv1d}
\lim_{h\to 0}\|A^h\|_\infty=0,\qquad \lim_{h\to 0} \|B^h\|_\infty=0.
\end{equation}
We need some preliminary lemmas.
\begin{Lemma}
Let  $\beta\ge 2$,
\begin{itemize}
\item[i)] For all $q, \tilde q\in \R^2$,
      \begin{eqnarray}
        \label{eq:29}
        \hskip-1cm    g(x,\tilde q)-g(x, q) -\nabla_q g(x, q)\cdot (\tilde q -q )
        &\ge & \frac 1 {\beta-1}  \max(|p|^{\beta-2},|\tilde p|^{\beta-2}) |p-\tilde p|^2
      \end{eqnarray}
where
      $ p=(q_1^-,q_2^+)$, $\tilde p=(\tilde q_1^-,\tilde q_2^+)$.
\item[ii)]
There exists a constant $C$ such that for all $q, \tilde q, r\in \R^2$ and $\eta>0$
\begin{equation}  \label{eq:32}
      \left|\left(\nabla_q g(x, \tilde q) -\nabla_q g(x, q)\right)\cdot r \right| \le         \max(|p|^{\beta-2},|\tilde p|^{\beta-2}) \left(\frac C \eta |p-\tilde p|^2     +  \eta |r|^2  \right),
    \end{equation}
where
      $ p=(q_1^-,q_2^+)$, $\tilde p=(\tilde q_1^-,\tilde q_2^+)$.
\end{itemize}
\end{Lemma}
For the proof of the previous lemma, we refer to \cite[Lemma 3.2]{acc}.
\begin{Lemma}
Let $(U^h,M^h,\Lambda^h)$ be a solution of  \eqref{MFGDs}-\eqref{transopm}  and $ (u^h, m^h, \l^h)$ a solution of \eqref{conv1a}-\eqref{conv1b} with $m^h\ge \d>0$ for   $h$ sufficiently small. Then
\begin{equation}\label{conv3}
    \lim_{h\to 0}\sum_{j\in J}\sum_{k=1}^{N^h_j-1}h_j\left| [D^+U]_{j,k} - [D^+u^h]_{j,k}\right|^\b=0
\end{equation}
\end{Lemma}
\begin{Proof}
By the identity \eqref{fundid} with $(U,M,\Lambda)=(U^h,M^h,\Lambda^h)$ and $(\bar U, \bar M, \bar \Lambda)=(u^h, m^h, \l^h)$ we get
\begin{multline*}
\cR^h(M^h,U^h,u^h)+\cR^h(m^h,u^h,  U^h)+(V_h[M^h]-V_h[m^h], M^h-m^h)_2
\\+(A^h,M^h-m^h)_2+(B^h,U^h-u^h)_2=0.\end{multline*}
By \eqref{stimasol} and the regularity of $u$ we get $ \lim_{h\to 0}  |(B^h,U^h-u^h)_2|=0$.
By $m^h, M^h\in \cK_h$ and the Cauchy-Schwarz
also get $ \lim_{h\to 0}  |(A^h,M^h-m^h)_2|=0$.
Hence by  \eqref{eq:29} we obtain
\begin{equation}
  \label{eq:46}
  \begin{array}[c]{rcl}
 \sum_{j\in J} \sum_{k=1}^{N_j^h} h_j\,m^h_{j,k}   \max\left\{|P^h_{j,k}|^{\beta-2}, |p^h_{j,k}|^{\beta-2} \right\} |P^h_{j,k}-p^h_{j,k}|^2 &=&o(1),\\
 \sum_{j\in J} \sum_{k=1}^{N_j^h}  h_j\,M^h_{j,k}  \max\left\{|P^h_{j,k}|^{\beta-2}, |p^h_{j,k}|^{\beta-2} \right\}  |P^h_{j,k}-p^h_{j,k}|^2 &=& o(1),
  \end{array}
\end{equation}
where $P^h_{j,k}=((D^+U^h)_{j,k}^-,(D^+U^h)_{j,k-1}^+)$ and $p^h_{j,k}=((D^+u^h)_{j,k}^-,(D^+u^h)_{j,k-1}^+)$. Since
$m^h$ is strictly positive, by the first equation in \eqref{eq:46} we get \eqref{conv3}.
\end{Proof}
\begin{Theorem}
Let $(u,m,\l)$ be the unique solution  of \eqref{MFGs} and
   $(U^h,M^h,\Lambda^h)$    the sequence of the solutions
of the scheme \eqref{MFGDs}-\eqref{transopm}. Then
\begin{equation}\label{convergence}
  \lim_{|h|\to 0} \|U^h-u\|_\infty+\|M^h-m\|_\infty+|\Lambda^h-\l|=0.
\end{equation}
\end{Theorem}
\begin{Proof}
We set $E^h=M^h-m^h$. Subtracting the equations satisfied by $M^h$ and $m^h$ and multiplying the resulting equations
for $h_jE^h_{j,k}$, we get
\begin{equation}\label{conv4}
\begin{split}
   \sum_{j\in J}  \sum_{k=1}^{N^h_j-1}\Big[-\nu_j (D^2_h E^h)_{j,k}+\cB^h(U^h,M^h)_{j,k}-\cB^h(U^h,m^h)_{j,k} \Big]h_jE^h_{j,k}=\\
-\sum_{j\in J}  \sum_{k=1}^{N^h_j-1}\big[  \cB^h(u^h,m^h)_{j,k} -\cB^h(U^h,m^h)_{j,k}+   B^h_{j,k}\big] h_jE^h_{j,k}.
\end{split}
\end{equation}
By the transition conditions for $M^h$ and $m^h$   (recall that  $M^h$ and, $m^h$ are continuous at the vertices)
\begin{equation}\label{conv5}
\begin{split}
&  \sum_{i\in I}\Big[\sum_{j\in \Inc^{+}_i}  E^h_{j,0}\big[\nu_j (D^+E^h)_{j,0}+E^h_{j,1}\frac{\pd g}{\pd q_2}(x_{j,1}, D_h U^h]_{j,1}\big]\\
&- \sum_{j\in \Inc^{-}_i}  E^h_{j,N^h_j}\big[\nu_j (D^+ E^h )_{j,N^h_j-1}+E^h_{j,N^h_j-1}\frac{\pd g}{\pd q_1}(x_{j,N^h_j-1}, [D_h U^h]_{j,N^h_j-1})\big]\Big]=\\
&\sum_{i\in I}\Big[\sum_{j\in \Inc^{+}_i}   E^h_{j,0}\big[m^h_{j,1}\big( \frac{\pd g}{\pd q_2}(x_{j,1}, [D_h u^h]_{j,1})-\frac{\pd g}{\pd q_2}(x_{j,1}, [D_h U^h]_{j,1})\big)+\frac{h_j}{2}B^h_{ j,0 }\big]
\\
&-\sum_{j\in \Inc^{-}_i}   E^h_{j,N^h_j}\big[ m^h_{j,N^h_j-1}\big(\frac{\pd g}{\pd q_1}(x_{j,N^h_j-1}, [D_h u^h]_{j,N^h_j-1})-\frac{\pd g}{\pd q_1}(x_{j,N^h_j-1}, [D_h U^h]_{j,N^h_j-1})\big)+\frac{h_j}{2}B^h_{ j,N^h_j }\big]\Big].
\end{split}
\end{equation}
 Arguing as in \eqref{ex2b}-\eqref{ex2c}, we have
\begin{align*}
 \sum_{j\in J}  \sum_{k=1}^{N^h_j-1}\nu_j (D^2_h E^h)_{j,k}E^h_{j,k}=  - \sum_{j\in J}  \sum_{k=1}^{N^h_j-2}\nu_j |(D^+ E^h)_{j,k}|^2
-\sum_{i\in I}\Big[\sum_{j\in \Inc^{+}_i} \nu_j | (D^+  E^h)_{j,0}|^2\\
+\sum_{j\in \Inc^{-}_i} \nu_j  |(D^+  E^h)_{j,N^h_j-1}|^2\Big]
+\sum_{i\in I}\Big[-\sum_{j\in \Inc^{+}_i}\frac{\nu_j}{h_j} E^h_{j,1} (D^+  E^h)_{j,0}
    +\sum_{j\in \Inc^{-}_i}\frac{\nu_j}{h_j} E^h_{j,N^h_j-1} (D^+  E^h)_{j,N^h_j-1}\Big].
\end{align*}
Moreover
\begin{align*}
\sum_{j\in J}\sum_{k=1}^{N^h_j-1} \Big[\cB^h(U^h,M^h)_{j,k}-\cB^h(U^h,m^h)_{j,k} \Big]E^h_{j,k}  =-\sum_{j\in J}\sum_{k=1}^{N^h_j-1}
E^h_{j,k}\,[D_h  E^h]_{j,k}\cdot \nabla_q g(x_{j,k},[D_hU^h]_{j,k})\\
 -\sum_{i\in I}\Big[\sum_{j\in \Inc^{+}_i}\frac{1}{h_j} E^h_{j,0}E^h_{j,1}\frac{\pd g}{\pd q_2}(x_{j,1},[D_hU^h]_{j,1}) -\sum_{j\in \Inc^{-}_i} \frac{1}{h_j} E^h_{j,N^h_j} E^h_{j,N^h_j-1}\frac{\pd g}{\pd q_1}(x_{j,N^h_j-1},[D_hU^h]_{j,N^h_j-1}) \Big].
\end{align*}
and
\begin{align*}
\sum_{j\in J}\sum_{k=1}^{N^h_j-1}&\Big[\cB^h(U^h,m^h)_{j,k}-\cB^h(u^h,m^h)_{j,k} \Big]E^h_{j,k} \\
 &=-\sum_{j\in J}\sum_{k=1}^{N^h_j-1} m^h_{j,k}\,[D_h  E^h]_{j,k}\cdot (\nabla_q g(x_{j,k},[D_hu^h]_{j,k})-\nabla_q g(x_{j,k},[D_hU^h]_{j,k})) \\
&-\sum_{i\in I}\Big[\sum_{j\in \Inc^{+}_i}\frac{1}{h_j} E^h_{j,0}m^h_{j,1}\left(\frac{\pd g}{\pd q_2}(x_{j,1},[D_hu^h]_{j,1})-
 \frac{\pd g}{\pd q_2}(x_{j,1},[D_hU^h]_{j,1})\right)
 \\
&-\sum_{j\in \Inc^{-}_i} \frac{1}{h_j} E^h_{j,N^h_j} m^h_{j,N^h_j-1}\left(\frac{\pd g}{\pd q_1}(x_{j,N^h_j-1},[D_hu^h]_{j,N^h_j-1})-\frac{\pd g}{\pd q_1}(x_{j,N^h_j-1},[D_hU^h]_{j,N^h_j-1})\right) \Big].
\end{align*}
Set
\[(D^h E^h, D^h E^h)_2 =\sum_{j\in J}  \sum_{k=1}^{N^h_j-2}  |(D^+ E^h)_{j,k}|^2
+\sum_{i\in I}\Big[\sum_{j\in \Inc^{+}_i}   | (D^+  E^h)_{j,0}|^2
+\sum_{j\in \Inc^{-}_i}   |(D^+  E^h)_{j,N^h_j-1}|^2\Big]\]
Replacing the previous   equalities in \eqref{conv4},  using \eqref{conv5}  and recalling  the estimate \eqref{erg}, we get
\begin{equation}\label{conv5b}
\begin{split}
&    (D^h E^h, D^h E^h)_2
 \le  - C\Big[
\sum_{j\in J}  \sum_{k=1}^{N^h_j-1} h_jE^h_{j,k}   A^h_{j,k}+h_j m^h_{j,k}\,[D_h  E^h]_{j,k}\cdot\big(\nabla_q g(x_{j,k},[D_hu^h]_{j,k})\\
 &-\nabla_q g(x_{j,k},[D_hU^h]_{j,k})\big) +\sum_{i\in I}\big(\sum_{j\in \Inc^{+}_i} h_j E^h_{j,0}
B^h_{ j,0 }
+\sum_{j\in \Inc^{-}_i} h_j E^h_{j,N^h_j}
B^h_{ j,N^h_j }\big)\Big]
\end{split}
\end{equation}
with $C$ independent of $h$. By \eqref{conv1d}, we have
\begin{equation}\label{conv6}
    \sum_{j\in J}  h_jE^h_{j,k} A^h_{j,k}+\sum_{i\in I}\Big[\sum_{j\in \Inc^{+}_i} h_j E^h_{j,0}B^h_{ j,0 }+\sum_{j\in \Inc^{-}_i} h_j E^h_{j,N^h_j}B^h_{ j,N^h_j }\Big]\le o(1) (E_h,E_h)_2.
\end{equation}
Set $ P^h=((D^+U^h_{j,k})^-,(D^+U^h_{j,k-1})^+)$, $p^h=((D^+u^h_{j,k})^-,(D^+u^h_{j,k-1})^+)$.
By  \eqref{stimasol}, \eqref{eq:32}  for any $\eta>0$
\begin{equation}\label{conv7}
\begin{split}
&\left|m^h_{j,k}\,[D_h  E^h]_{j,k}\cdot (\nabla_q g(x_{j,k},[D_hu^h]_{j,k})-\nabla_q g(x_{j,k},[D_hU^h]_{j,k}))\right|\\
& \le   m^h_{j,k}\max(|P^h_{j,k}|^{\beta-2},|  p^h_{j,k}|^{\beta-2}) \left(\frac C \eta |P^h_{j,k}- p^h_{j,k}|^2     +  \eta |D^hE^h_{j,k}|^2  \right)\\
&\le    m^h_{j,k} \left(\frac C \eta |[D^h U^h]_{j,k}-[D^h u^h]_{j,k} |^\beta     +  \eta |D^hE^h_{j,k}|^2  \right).
\end{split}
\end{equation}
Plugging the estimates \eqref{conv3}, \eqref{conv6}, \eqref{conv7} in \eqref{conv5b}, we finally get
\[ (E^h, E^h)_2+(D^h E^h, D^h E^h)_2=o(1) \quad\text{for $|h|\to 0$}.\]
Hence we get  the convergence of $M^h$ to $m$ in $H^1(\G)$ and uniform.
By the convergence of $M^h$ to $m$ and \textbf{(V4)}, we get $\lim_{h\to 0}\|V_h[M^h]-V_h[m^h]\|_\infty=0$. Hence
$U^h$ and $u^h$ are solution of \eqref{HJBs} with $\l=\Lambda^h$, $V_{j,k}=V_h[M^h]_{j,k}$ and respectively
$\l=\l^h$, $V_{j,k}=V_h[  M^h]_{j,k}+o(1)$. By a comparison principle  for \eqref{HJBs}, we get
$|\Lambda^h-\l^h|\le o(1)$ and therefore
\begin{equation}\label{conv9}
\lim_{|h|\to 0}|\l-\Lambda^h|=0.
\end{equation}
Let $\bar u^h$ be the continuous piecewise linear reconstruction of $U^h$ on $\G$. By \eqref{stimasol}, $\bar u^h\to \bar u$
uniformly as $|h|\to 0$, up to a subsequence. By \eqref{conv3} and \eqref{conv9}, $\bar u$ is a weak solution to \eqref{MFGs}. Therefore
by the uniqueness of the solution to \eqref{MFGs}, we get the convergence of $U^h$ to $u$ in $H^1(\G)$ and uniform.
\end{Proof}

\section{Numerical implementation and experiments}\label{sec4}
This section is devoted to the   implementation and test of a numerical solver for the stationary MFG system
\eqref{MFGDs}-\eqref{transopm}.
In \cite{acd},  the stationary MFG system on the torus
is solved  via the so called {\em forward-forward long time approximation}:
for a given approximation step $h$, the approximate solution  $(U_h,M_h,\Lambda_h)$   is obtained as the limit of
$(U^n_h, M^n_h, U^n_h/n\Delta t) $
for $n\to \infty$, where $(U^n_h,M^n_h)$ is computed  via discretization of the corresponding  evolutive MFG system, implicit or explicit  in time, up to time $T=n\Delta t$.\par
Here we propose a new approach which allows  to compute the solution of the stationary MFG system {\em directly}, avoiding    {\em long time} or {\em small delta} approximations. We   collect all the unknowns $(U,M,\Lambda)$ in a single vector $X$ of length $2N^h+1$ (with $N^h$ given by \eqref{hi}) and we recast   the $2N^h+2$ equations of the
stationary MFG system as functions of $X$. Hence  we get a  nonlinear map
$\mathcal F:\mathbb{R}^{2N^h+1}\to\mathbb{R}^{2N^h+2}$ defined by
$$\mathcal F(X)=
\left\{
\begin{array}{ll}
    -\nu_j (D^2_h U)_{j,k} + g(x_{j,k}, [D_h U ]_{j,k} )+\Lambda - (V_h[M])_{j,k} & k=1,\dots, N^h_j-1,\, j\in J\\[4pt]
     \nu_j(D^2_h M)_{j,k}+\cB^h(U,M)_{j,k} & k=1,\dots, N^h_j-1,\, j\in J\\[4pt]
     \cS^h(U,V_h[M]-\Lambda)_i  &i\in I\\[4pt]
     \cT^h(M,U)_i &i\in I\\[4pt]
     (M,1)_2-1 \\ (U,1)_2
\end{array}
\right.
$$
and we look for $X^\star\in\mathbb{R}^{2N^h+1}$ such that
\begin{equation}\label{S1}
\mathcal F(X^\star)=0\in\mathbb{R}^{2N^h+2}.
\end{equation}
By Theorem \ref{th_existence} and  Theorem \ref{th_uniq} there exists a unique solution to \eqref{S1}, but the system is formally {\em overdetermined}, having $2N^h+2$ equations in $2N^h+1$ unknowns.
{\bf This terminology normally applies to linear systems, but is commonly adopted also in the nonlinear case with a slight abuse of notation. }
Indeed, the solution is meant in the following {\em nonlinear-least-squares} sense:
$$X^\star=\arg\min_{X}\frac12\|\mathcal F(X)\|_2^2\,.$$
To solve the above optimization problem, we employ the {\em Gauss-Newton} method, that we briefly
recall here for completeness. We first denote the {\em residual function} by
$$r(X)=\frac12\|\mathcal F(X)\|_2^2=\frac12 \mathcal F(X)^T\mathcal F(X)$$ and we consider the standard Newton method for approximating a
critical point of $r$:
$$\cH_r(X^k)\delta_X=-\nabla r(X^k)\,,\qquad X^{k+1}=X^k+\delta_X\,,\qquad k\ge 0\,,$$
where the gradient $\nabla r$ and the Hessian $\cH_r$ are given by
$$\nabla r(X)=J_\mathcal F(X)^T \mathcal F(X)\,,\qquad \cH_r(X)=J_\mathcal F(X)^T J_\mathcal F(X)+\sum_{i=1}^{2N^h+2}\frac{\partial^2 \mathcal F_i}{\partial^2 X}(X)
\mathcal F_i(X)\,,$$
with
$$\left(J_\mathcal F(X)\right)_{i,j}=\frac{\partial \mathcal F_i}{\partial X_j}(X)\,,\qquad \left(\frac{\partial^2 \mathcal F_i}{\partial^2 X}(X)\right)_{k,\ell}=\frac{\partial^2 \mathcal F_i}{\partial X_k \partial X_\ell}(X)\,.$$
Since we expect the residuals $\mathcal F_i(X^k)$ to be small for $X^k$ close enough to $X^\star$, it is reasonable
to neglect the second derivatives in $\cH_r$, using the approximation $\cH_r(X)\simeq J_\mathcal F(X)^T J_\mathcal F(X)$.
This yields the Gauss-Newton method:
$$J_\mathcal F(X^k)^T J_\mathcal F(X^k)\delta_X=-J_\mathcal F(X^k)^T\mathcal F(X^k)\,,\qquad X^{k+1}=X^{k}+\delta_X\,,\qquad k\ge 0\,,$$
where the Jacobian $J_\mathcal F$ is well defined assuming that the numerical Hamiltonian $g$
is of class $C^2$ in the gradient variable and that the operator $V_h$ is of class $C^1$.
Despite this method allows to employ only first order derivatives of $\mathcal F$, it is still
not efficient from a numerical point of view. Indeed, once $J_\mathcal F$ at $X^k$ is computed,
we also need to assemble the right hand side $J_\mathcal F^T \mathcal F$ and the matrix $J_{\mathcal F}^T J_\mathcal F$, typically squaring
the condition number of the system.
This can be avoided by simply realizing that the $k$-th iteration of the Gauss-Newton method is just
the normal equation for the following {\em linear-least-squares} problem
\begin{equation}\label{leastsquares}
\min_{\delta_X}\frac12\|J_\mathcal F(X^k)\delta_X+\mathcal F(X^k)\|_2^2\,,
\end{equation}
which is in turn easily and efficiently solved by means of the {\em $QR$ factorization} of $J_\mathcal F$. Indeed,
let $m=2N^h+2$, $n=2N^h+1$ and suppose that $J_\mathcal F(X^k)=QR$, where $Q$ is a $m\times m$
orthogonal matrix (i.e. $Q^{-1}=Q^T$) and $R$ is a $m\times n$ matrix of the form
$R=\left(\begin{array}{c} R_1 \\ 0 \end{array}\right)$, with $R_1$ of size $n\times n$ and upper triangular.
Writing $Q=\left( Q_1\quad Q_2\right)$ with $Q_1$ of size $m\times n$ and $Q_2$ of size $m\times (m-n)$, we get
$$
\|J_\mathcal F(X^k)\delta_X+\mathcal F(X^k)\|_2^2=\|Q^T\left(J_\mathcal F(X^k)\delta_X+\mathcal F(X^k)\right)\|_2^2
=\|Q^T QR\delta_X+Q^T\mathcal F(X^k)\|_2^2=
$$
$$
=\left\| \left(\begin{array}{c} R_1 \delta_X \\ 0 \end{array}\right)+
 \left(\begin{array}{c} Q_1^T \mathcal F(X^k) \\ Q_2^T \mathcal F(X^k)\end{array}\right)
\right \|_2^2=
\|R_1\delta_X+Q_1^T\mathcal F(X^k)\|_2^2+\|Q_2^T \mathcal F(X^k)\|_2^2
$$
which is finally minimized by getting rid of the first of the two latter terms,
i.e. solving the square triangular $n\times n$ linear system $R_1\delta_X=-Q_1^T\mathcal F(X^k)$
via back substitution. \par

Summarizing, we propose the following simple algorithm for the stationary MFG system:\\[4pt]
\textsc{
Given a guess $X=(U^0,M^0,\Lambda^0)$, a tolerance $\varepsilon>0$ and a dumping parameter $0<\alpha\leq 1$,\\
repeat
\begin{itemize}
 \item[$\bullet$] Assemble $\mathcal F(X)$ and $J_\mathcal F(X)$
 \item[$\bullet$] Solve the overdetermined linear system $J_\mathcal F(X)\delta_X=-\mathcal F(X)$
 in the least-squares sense \eqref{leastsquares},
 using the $QR$ factorization of $J_\mathcal F(X)$
 \item[$\bullet$] Update $X\,\leftarrow\,X+\alpha\delta_X$
\end{itemize}
until $\|\delta_X\|_2<\varepsilon$}\\

\noindent The algorithm is implemented in    \textsc{C}-language and  employs the library {\em SuiteSparseQR}
\cite{SPQR}, which is designed to efficiently compute the $QR$ factorization and the least-square solution to
very large and sparse linear systems.\par

\noindent Some remarks are in order:\\
  1)
We always initialize the method by setting $U^0\equiv 0$, $\Lambda^0=0$ and $M^0\equiv 1/L$, where $L=\sum_{j\in J}l_j$ is the total length of the network. In general there is no guarantee that the algorithm computes a
 minimum of  \eqref{leastsquares} with zero residual, i.e. a solution of the stationary MFG system. Nevertheless,
in all the tests performed, our algorithm seems to converge  to a zero residual minimum
independently on the initial guess. \par
\noindent 2) As for the standard Newton method, it is known that also the Gauss-Newton method may not converge
if the dumping parameter is set to $\alpha=1$. A fine tuning of $\alpha$ can be accomplished via some moderate
time consuming line search technique, but for our purposes we simply checked that the fixed value $\alpha=0.9$ is sufficient
in all the considered examples.\par
\noindent 3) We never impose  the constraint $M\ge 0$ in the computation. Surprisingly, our unconstrained optimization
algorithm converges to a solution of the MFG system with non negative mass.
We extensively checked this feature, also in the case of negative or changing sign initial guesses. Even if
the mass can be negative in some intermediate iterations of the Gauss-Newton method,
we always end up with a non negative mass in all the considered examples.\par
\noindent 4) Our technique can be successfully applied also in the homogenization of Hamilton-Jacobi equations, e.g.
for computing the effective Hamiltonian for some cell problems. Our preliminary tests using the nonlinear-least-squares approach are very promising, both in terms of accuracy and computational costs. \par
The previous points, in particular the convergence of the method, are still under investigation and will be
addressed in a future work (see \cite{cc}).\\[4pt]

\begin{figure}[h!]
 \begin{center}
 \includegraphics[width=.7\textwidth]{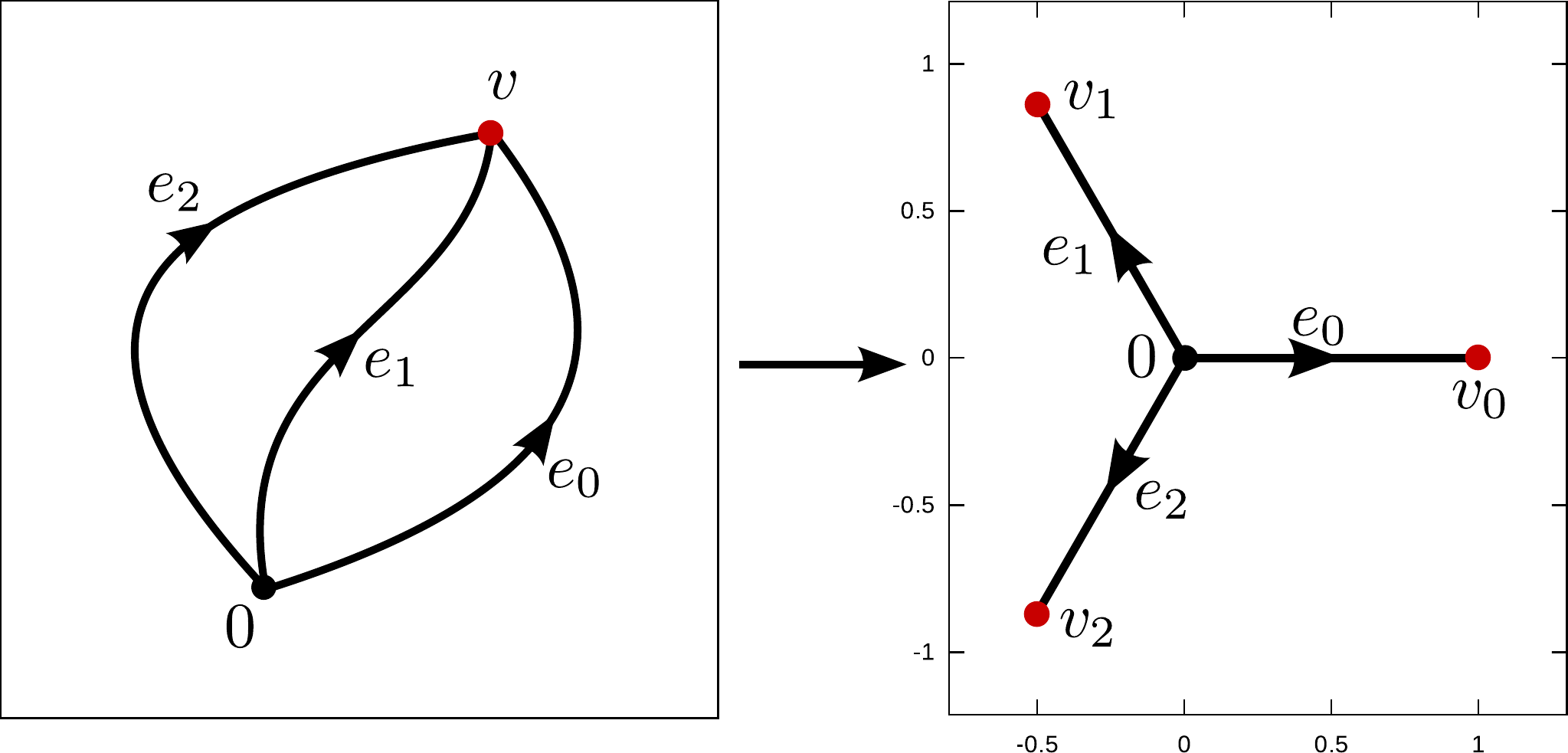}
 $$(a)\qquad\hskip1.5cm\qquad\qquad\qquad\qquad\qquad(b)$$
 \end{center}
 \caption{a network with 2 vertices and 3 edges (a) is mapped in an equivalent network with boundary vertices identified (b).}\label{net3edges}
 \end{figure}

We now set up the data for the numerical experiments. We consider a simple network in the plane
with $2$ vertices and $3$ edges of unit length, as in Figure \ref{net3edges}a.\textbf{ For computational purposes the network is mapped in
a topologically equivalent network, in which one vertex is located at the origin and the edges are delimited by the $3$rd roots of unity
($v_j=(\cos(2\pi j/3),\sin(2\pi j/3))$, for $j=0,1,2$), as in Figure \ref{net3edges}b.
Note that the boundary vertices, i.e. the vertices with a single incident edge, are identified and correspond to a single vertex on the network.}
\\
We assume that the numerical Hamiltonian has the form \eqref{hamnumconv}, with $\beta=2$ and $f(x)$ is such that,
for $j=0,1,2$ and $x\in e_j$,
$$f(x)=f_j(x):=s_j\left(1+\cos\left(2\pi(t+1/2)\right)\right)\,,\qquad x=tv_j\,,\quad t\in[0,1]\,,$$
where $s_j\in\{0,1\}$ is a switch parameter to activate/deactivate the corresponding cost on the edge $e_j$.
If not differently specified, we discretize each edge by $N_j=250$ nodes, so that the resulting nonlinear system has dimension
$1502\times 1501$, and we choose a tolerance $\varepsilon=10^{-4}$ for the stopping criterion of the algorithm.
We finally assume a uniform diffusion on the whole network, i.e. $\nu_j\equiv\nu$ for $j=0,1,2$ and $\nu>0$.\\
{\bf Here we are mainly interested in the qualitative behavior of the computed solutions,
and we postpone at the end of the section some experimental analysis on the performance of the algorithm.} Nevertheless, we remark that in all
the following tests, the proposed method converges in about $10$ iterations and the computational time
is of the order of few seconds, even for larger grids.\\ 

\noindent{\bf All the tests were performed on a Lenovo Ultrabook X1 Carbon, using 1 CPU Intel Quad-Core i5-4300U 1.90Ghz with 8 Gb Ram, 
running under the Linux Slackware 14.1 operating system.} \\

\begin{figure}[h!]
 \begin{center}
 \begin{tabular}{ccc}
 \footnotesize$\,\,\,\,\,\,\,\min=0.039\,\,\,\,\,\,\,\,\,\,\max=0.778$ &
 \footnotesize$\,\,\,\,\,\,\,\min=5\!\times\!\!10^{-4}\,\,\,\,\,\max=1.017$ &
 \footnotesize$\,\,\,\,\,\,\,\min=0.053\,\,\,\,\,\,\,\,\,\,\,\max=1.328$\\
 \quad\includegraphics[width=.275\textwidth]{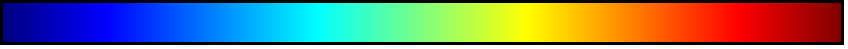} &
 \quad\includegraphics[width=.275\textwidth]{mfg-colorbar} &
 \quad\includegraphics[width=.275\textwidth]{mfg-colorbar} \\
\includegraphics[width=.3\textwidth]{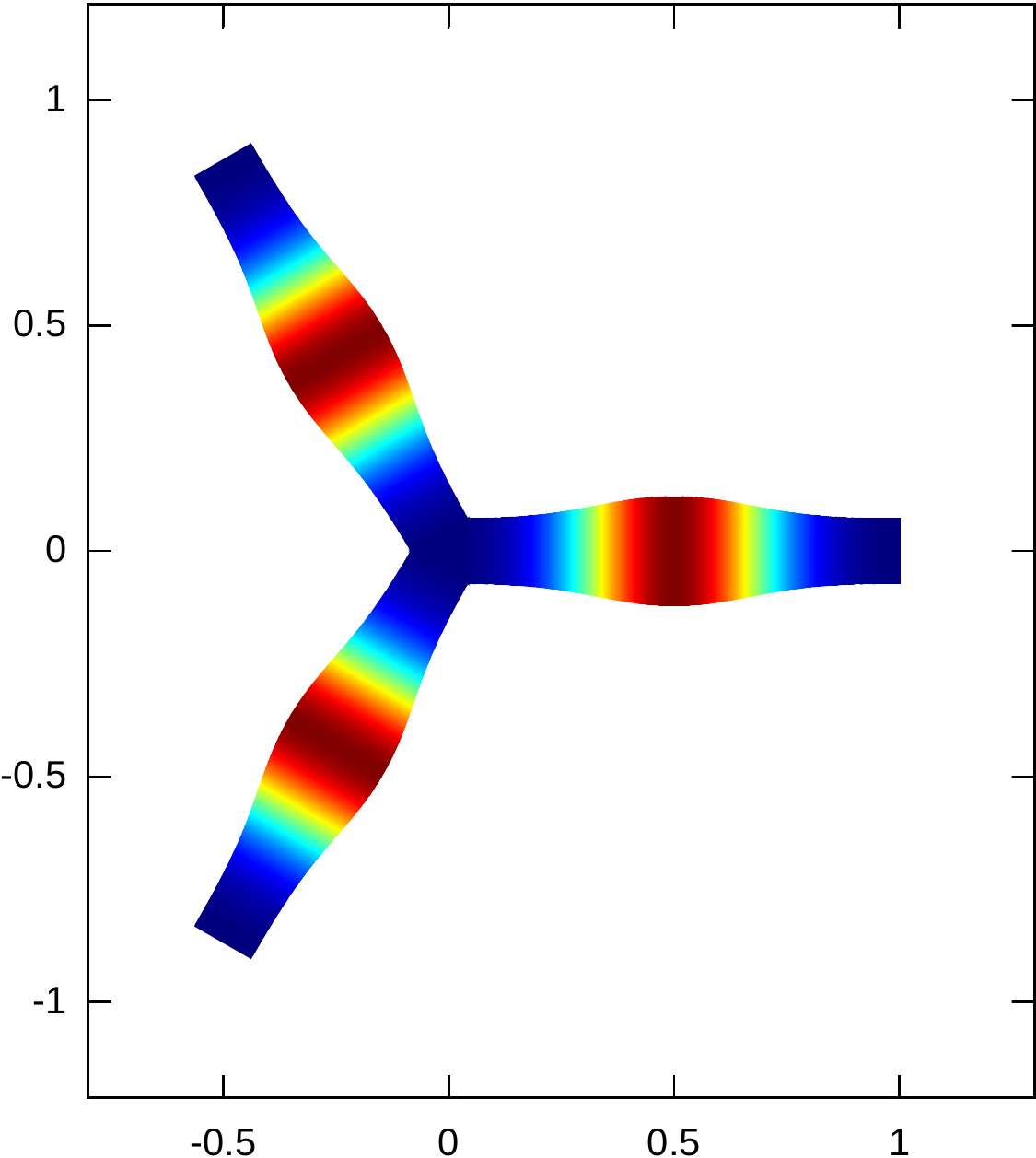} &
\includegraphics[width=.3\textwidth]{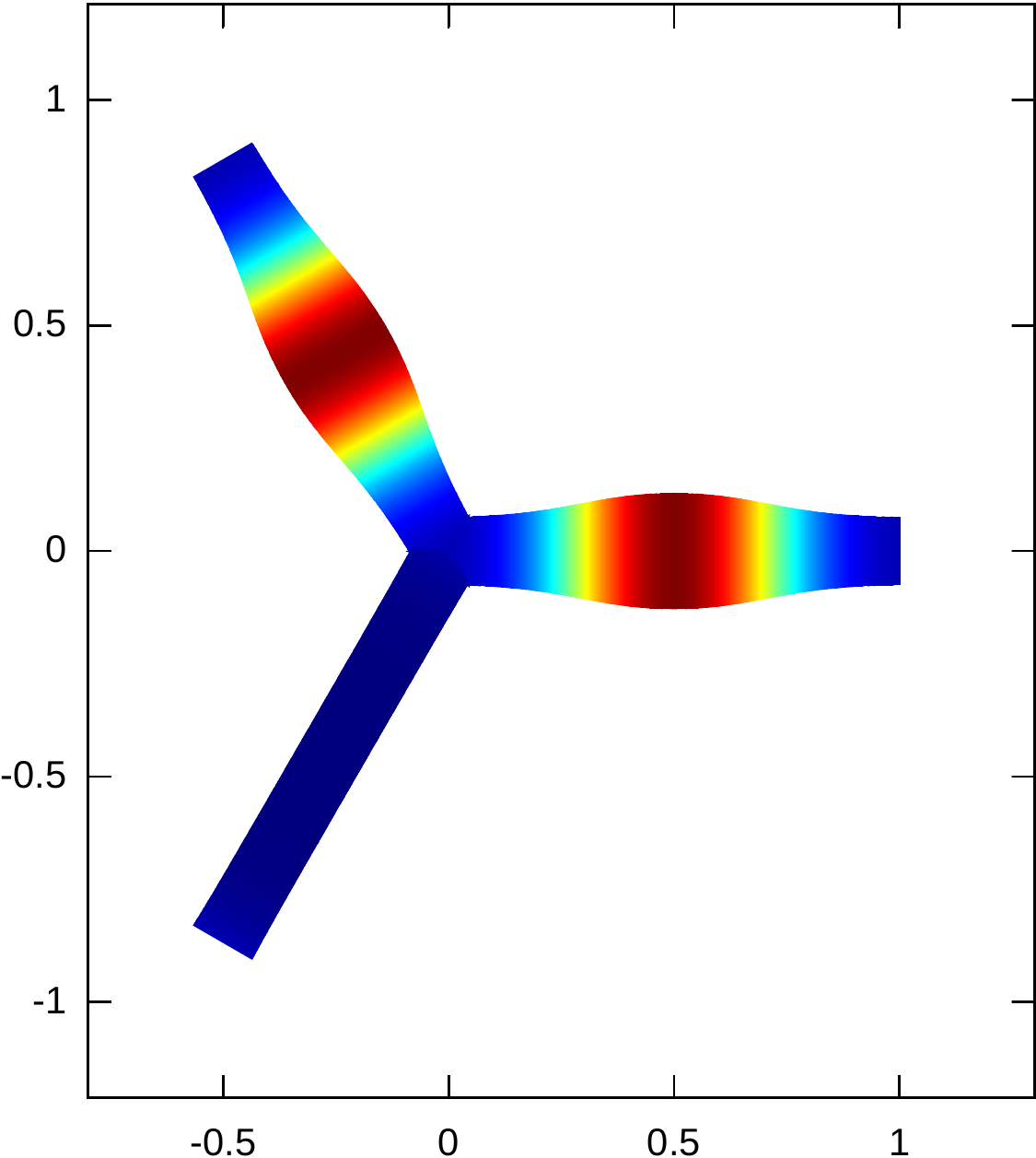} &
\includegraphics[width=.3\textwidth]{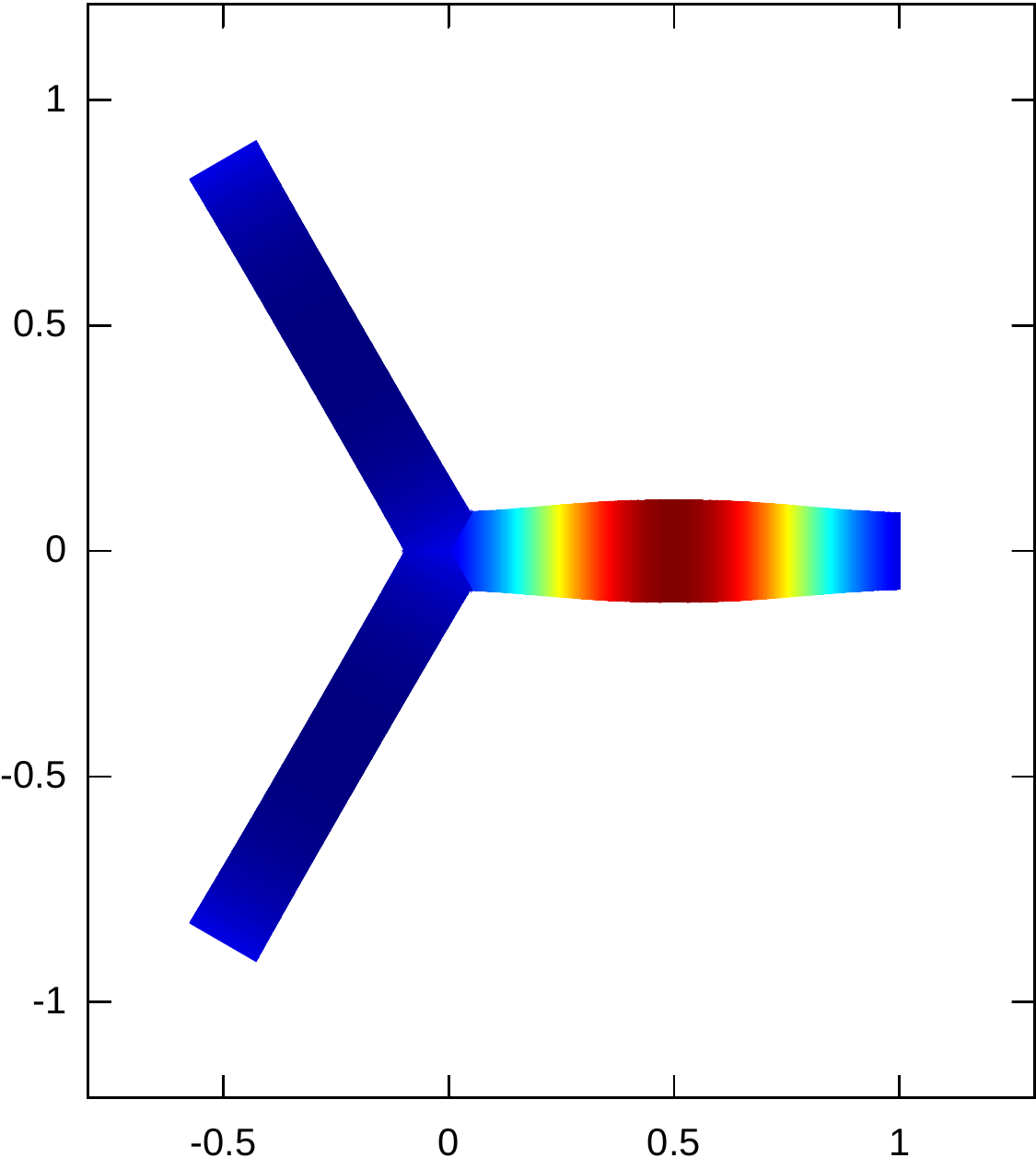} \\
\\
\includegraphics[width=.3\textwidth]{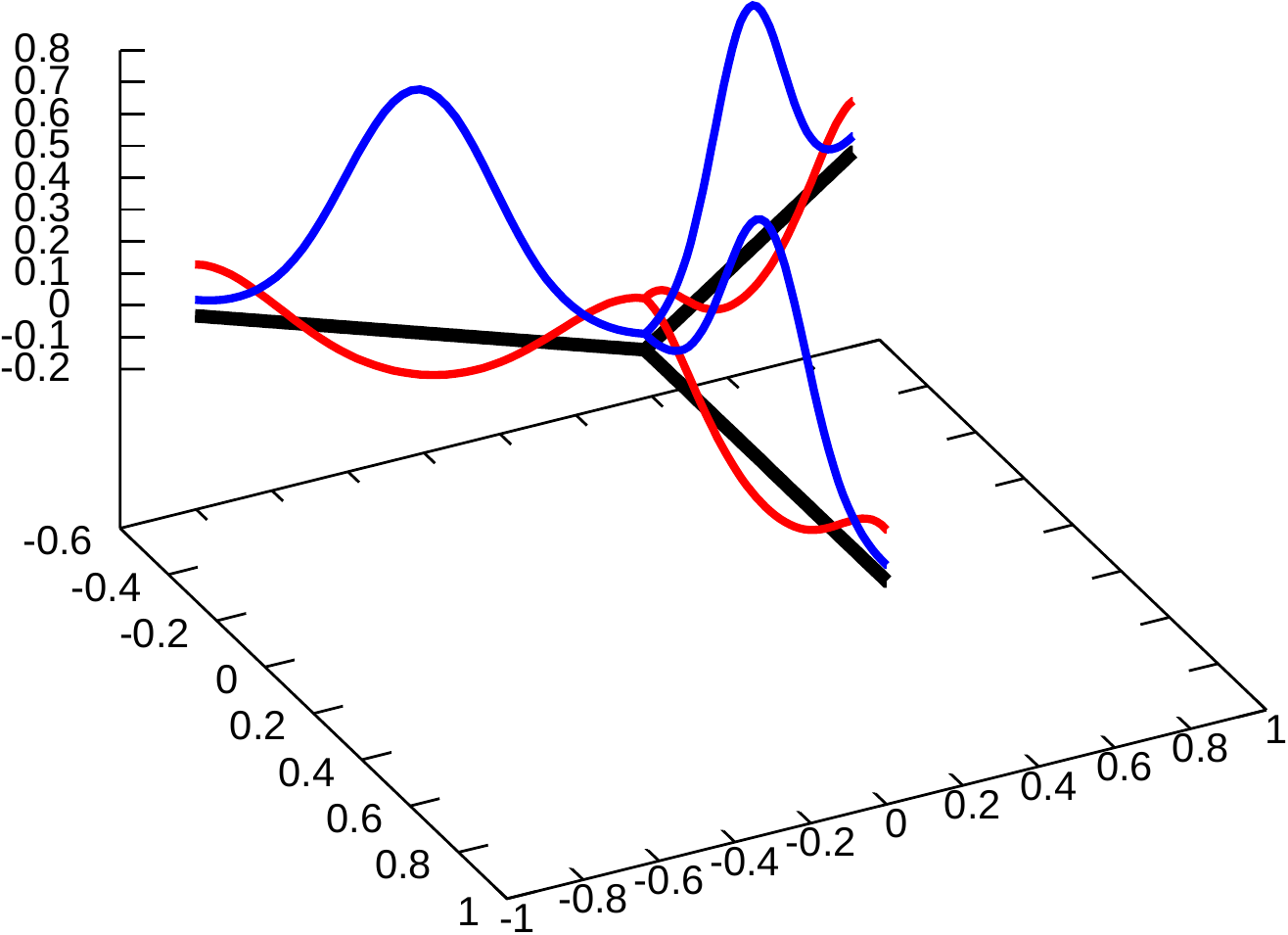} &
\includegraphics[width=.3\textwidth]{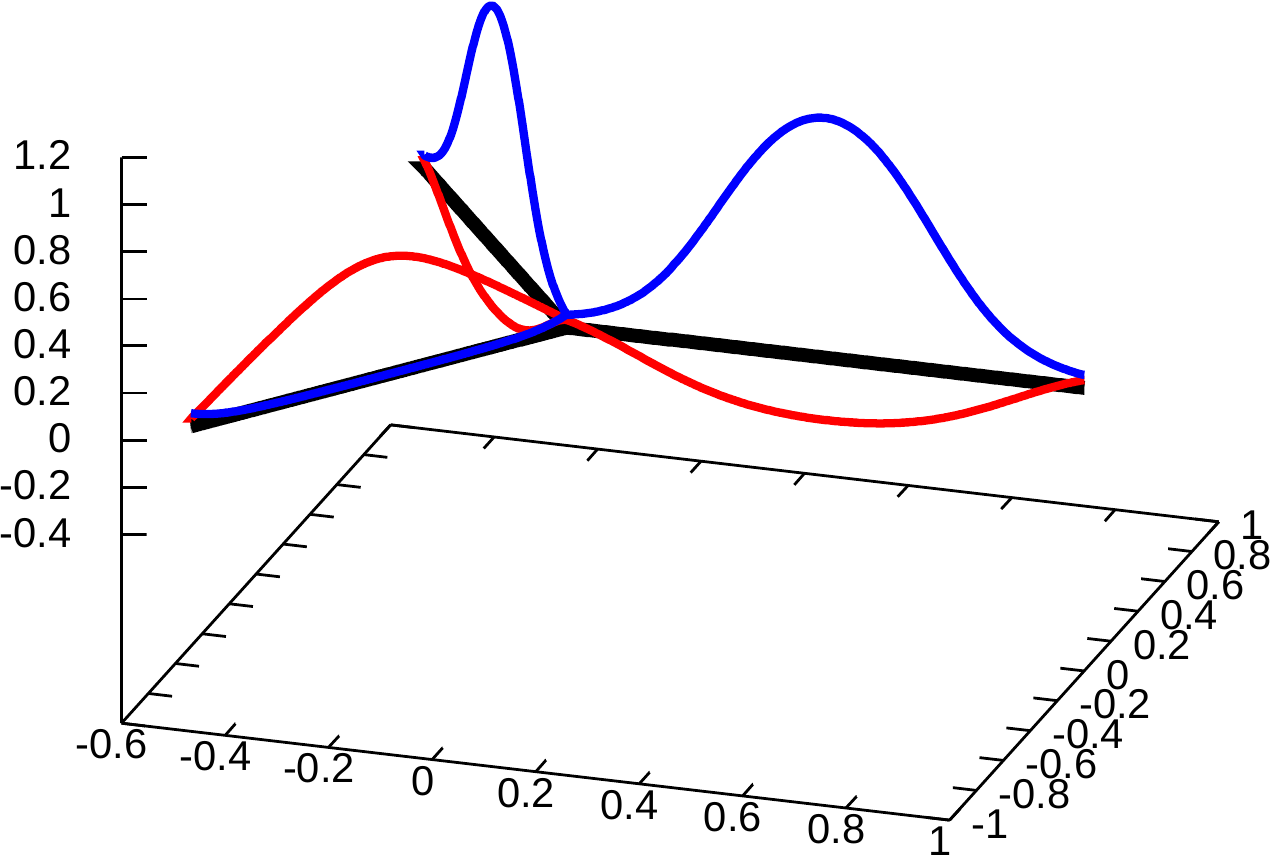} &
\includegraphics[width=.3\textwidth]{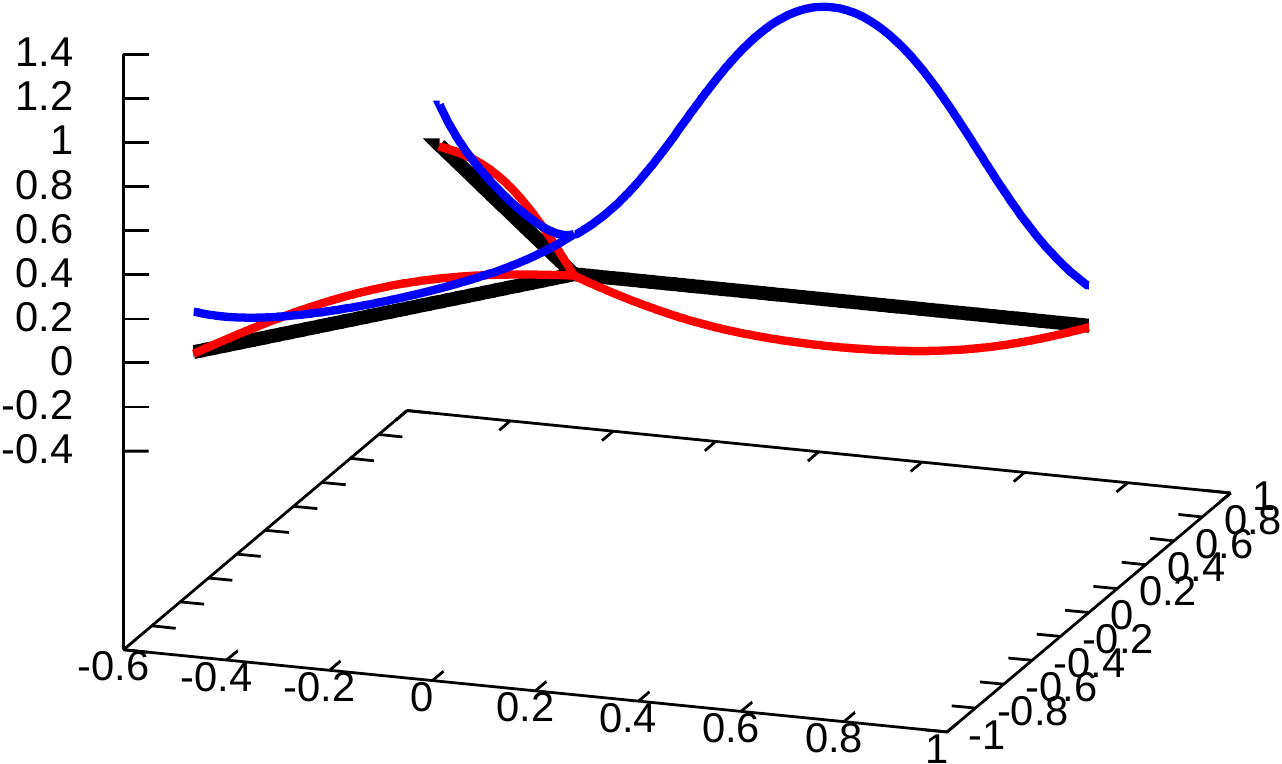} \\
(a)&(b)&(c)
\end{tabular}
\end{center}
\caption{the case $V_h[m]=m^2$ and $\nu=0.1$, the cost $f$  is active on (a) three edges, (b) two edges, (c) one edge.}\label{test1}
\end{figure}

\noindent{\bf Test 1.} We consider a local operator of the form $V_h[m]=m^2$,
and we choose a diffusion coefficient $\nu=0.1$. Figure \ref{test1} shows the results corresponding to the activation of the cost
$f$ on three, two or one edge, namely for $(s_0,s_1,s_2)=(1,1,1)$, $(s_0,s_1,s_2)=(1,1,0)$
and $(s_0,s_1,s_2)=(1,0,0)$ respectively. 
\\
In the top panels we represent the mass $M$ using a color-map in which
the blue and the red correspond respectively to the  minimum and maximum values. Moreover, we represent the network as
a fatten tube, whose cross sections have a size proportional to $M$ at the corresponding points.
In the bottom panels we represent the network (in black) and both the mass $M$ (in blue) and the corresponding value function $U$ (in red).
Since   $V$ is increasing, it penalizes concentration of the mass. The cost  $f$ has, if $s_j\neq 0$,  a maximum in the center of the edge $e_j$. Hence, if $\nu$ is not to small, the agents should be well distributed on the network with a maximum of $m$ around the minima of the value functions  $u$, i.e. in the center of the edges where the cost is active. In fact we observe    this  behavior in all the three examples.\\

\begin{figure}[h!]
 \begin{center}
 \begin{tabular}{ccc}
  \footnotesize$\,\,\,\,\,\,\,\min=0\,\,\,\,\,\,\,\,\,\,\,\,\,\,\,\,\,\,\,\,\,\max=0.939$ &
 \footnotesize$\,\,\,\,\,\,\,\min=0\,\,\,\,\,\,\,\,\,\,\,\,\,\,\,\,\,\,\,\,\,\max=1.129$ &
 \footnotesize$\,\,\,\,\,\,\,\min=0\,\,\,\,\,\,\,\,\,\,\,\,\,\,\,\,\,\,\,\,\,\max=1.915$\\
 \quad\includegraphics[width=.275\textwidth]{mfg-colorbar} &
 \quad\includegraphics[width=.275\textwidth]{mfg-colorbar} &
 \quad\includegraphics[width=.275\textwidth]{mfg-colorbar} \\
\includegraphics[width=.3\textwidth]{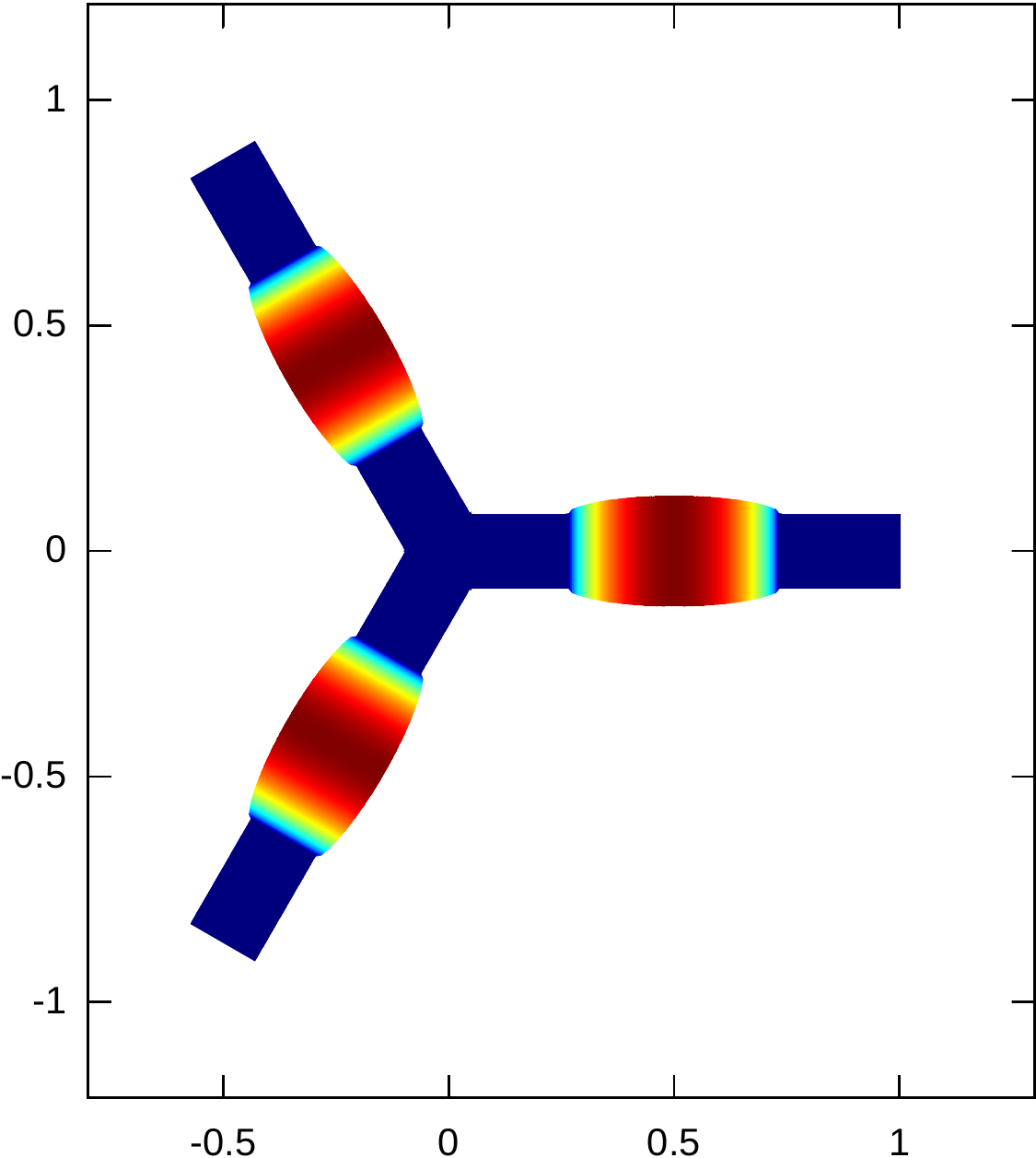} &
\includegraphics[width=.3\textwidth]{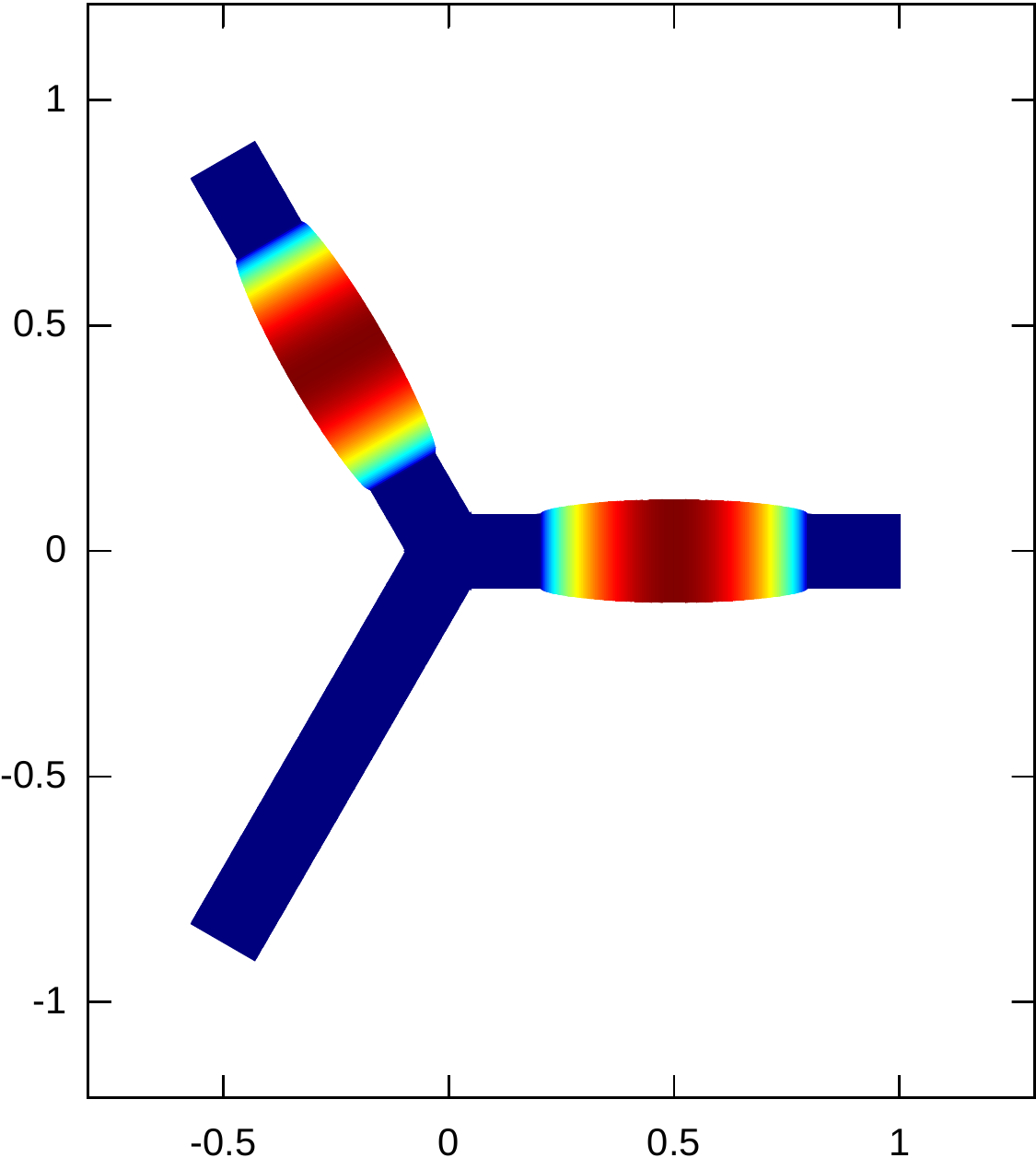} &
\includegraphics[width=.3\textwidth]{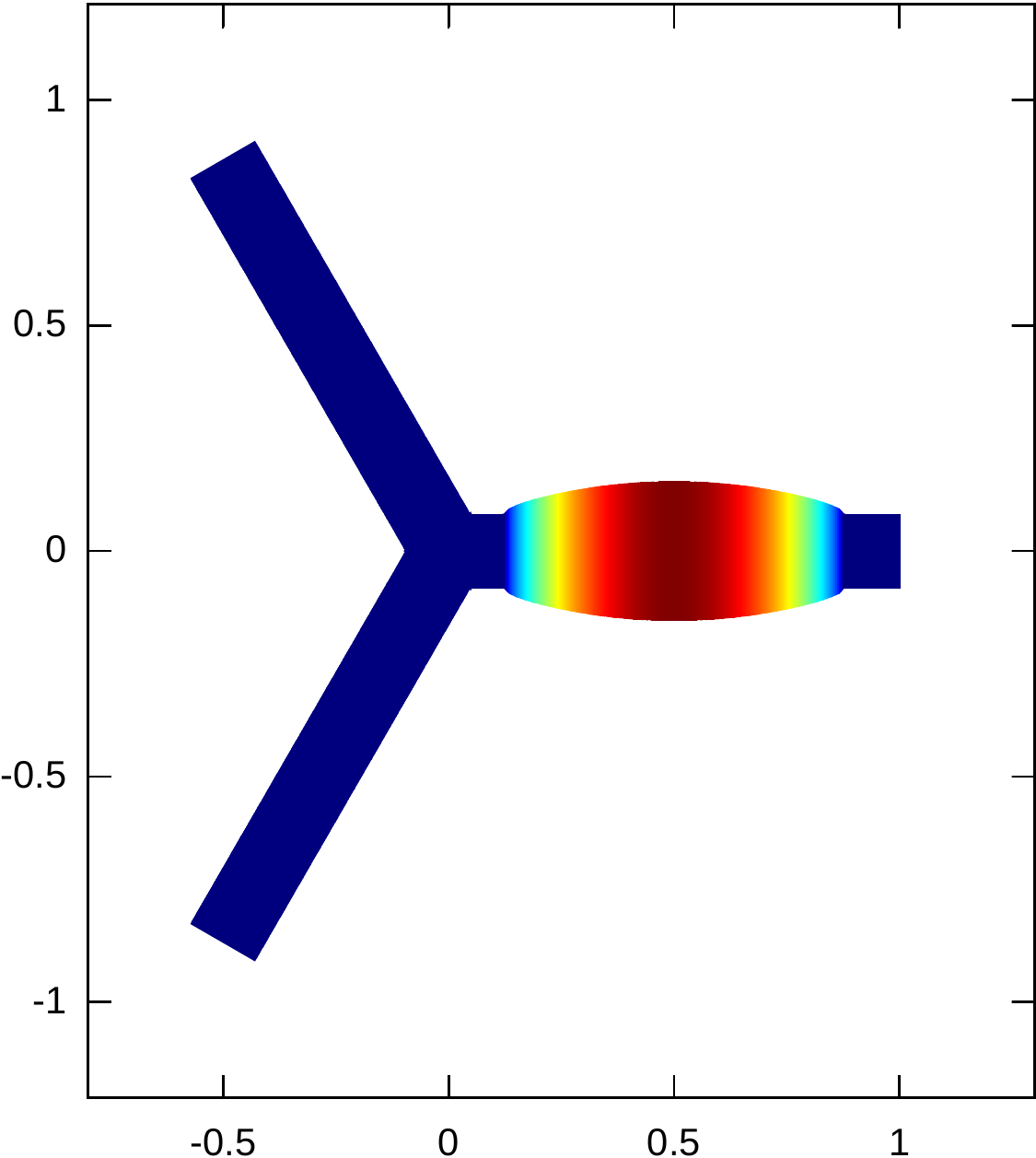} \\
\\
\includegraphics[width=.3\textwidth]{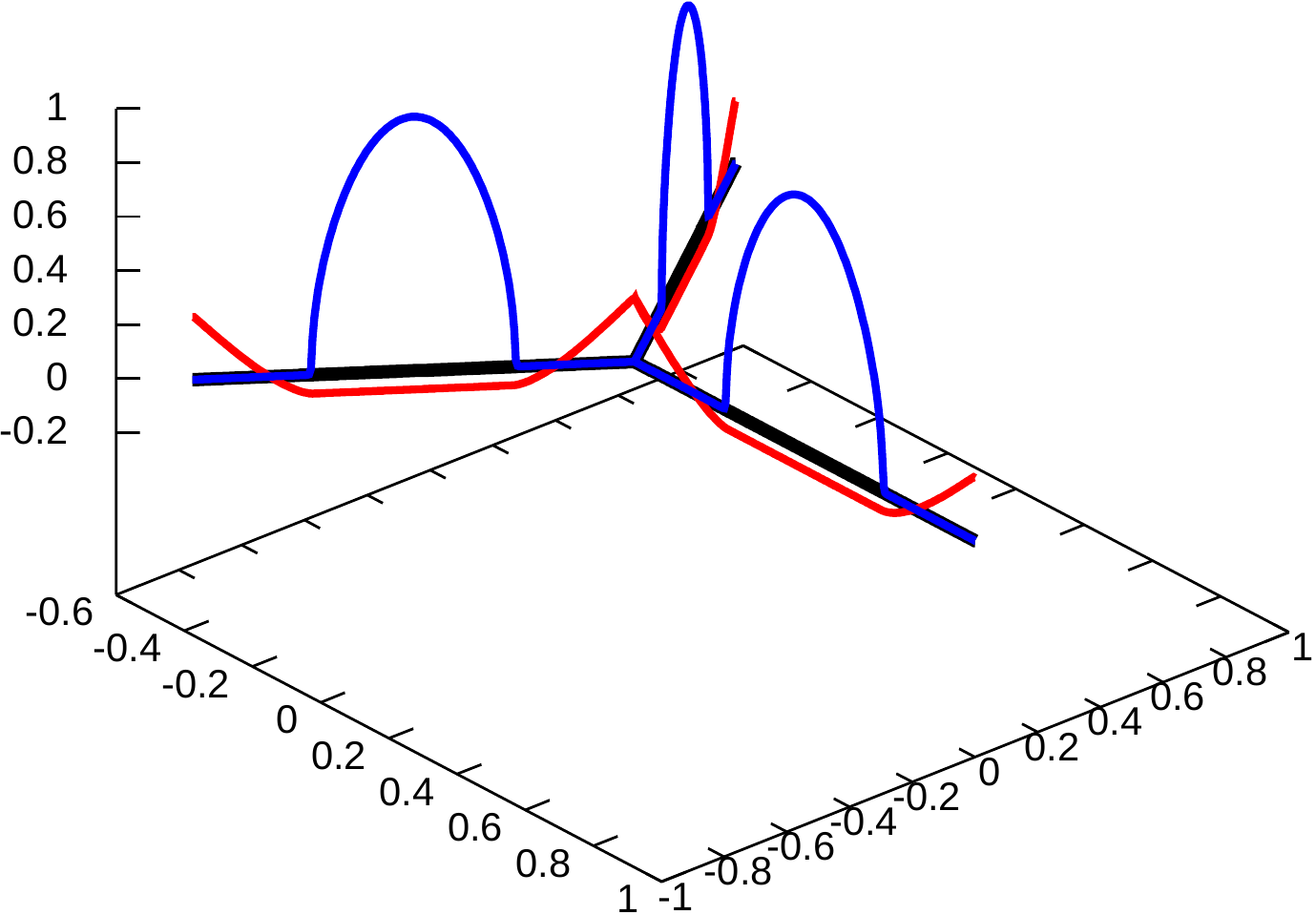} &
\includegraphics[width=.3\textwidth]{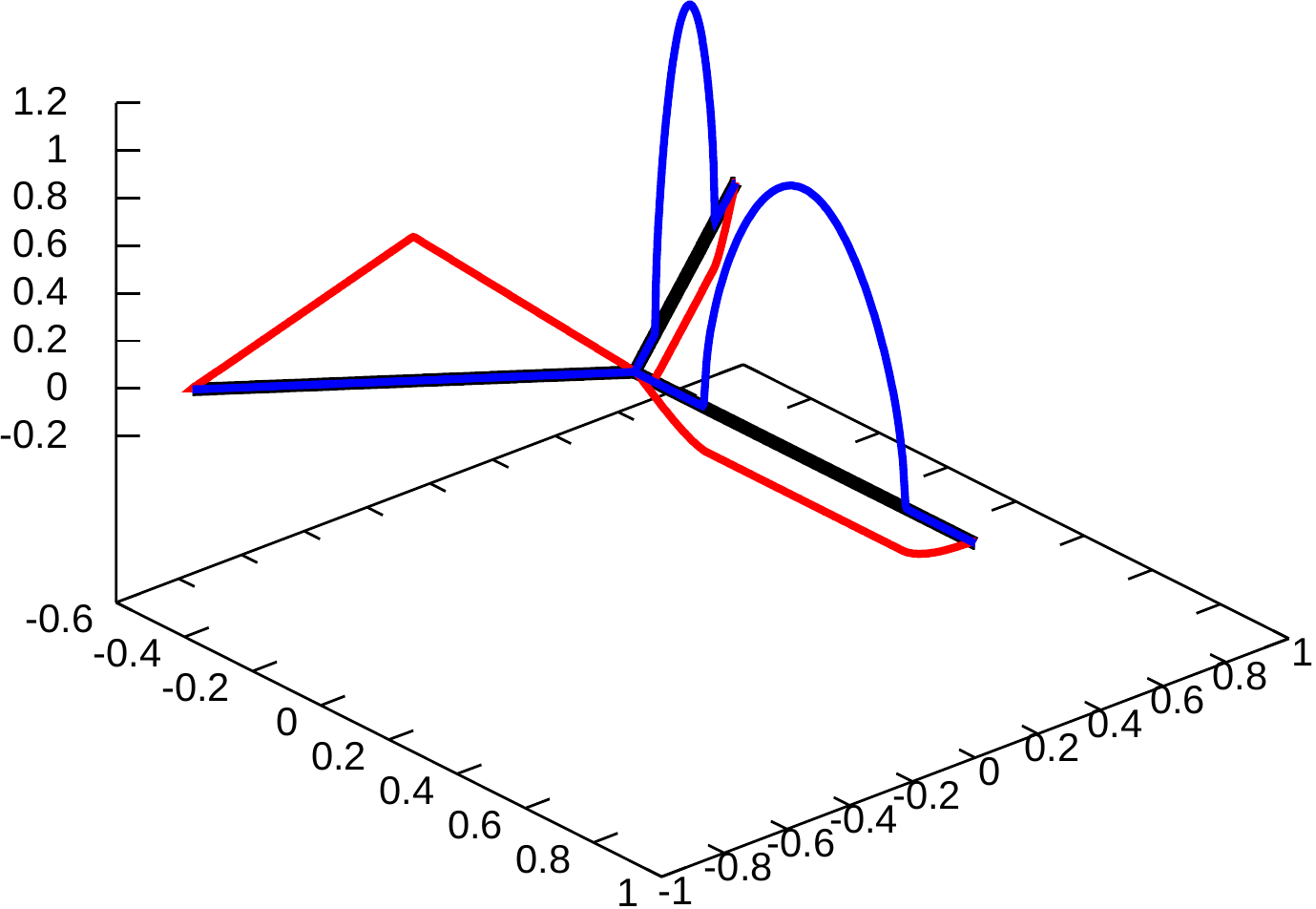} &
\includegraphics[width=.3\textwidth]{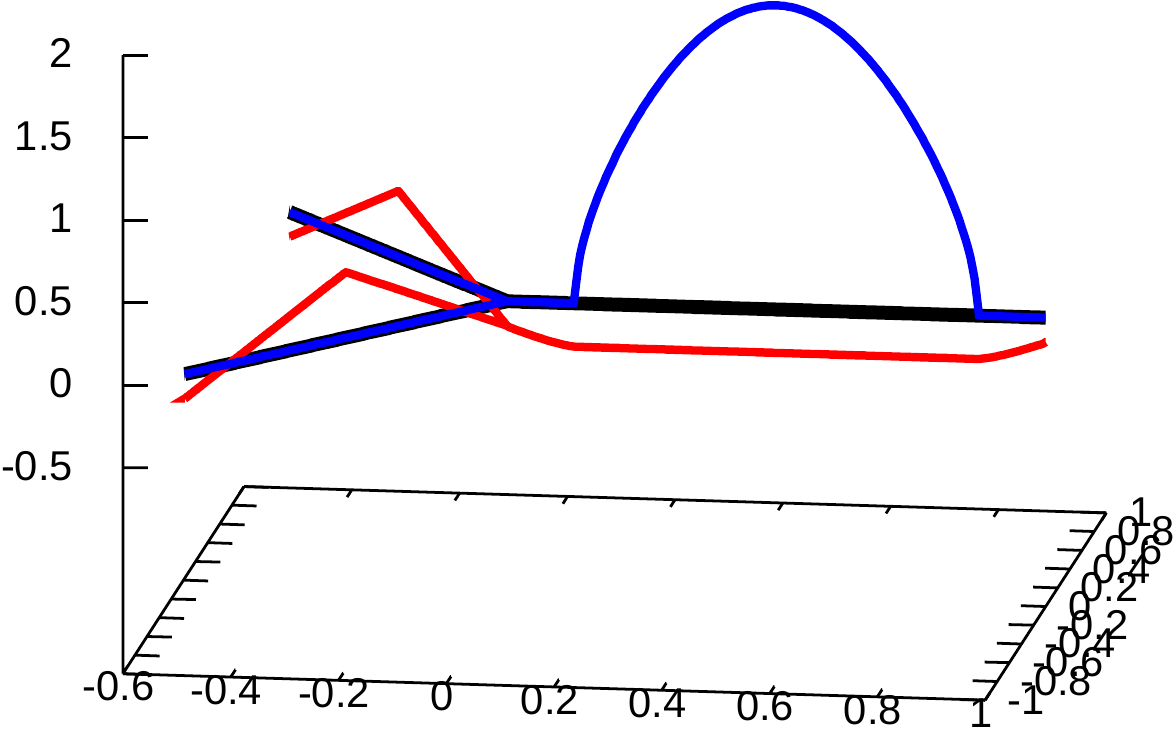} \\
(a)&(b)&(c)
\end{tabular}
\end{center}
\caption{the case $V_h[m]=m^2$ and $\nu=10^{-4}$, the cost $f$ is active on (a) three edges, (b) two edges, (c)
one edge.}\label{test2}
\end{figure}

\noindent{\bf Test 2.} We are interested in the behavior of the solution  as $\nu\to 0$, hence
we choose the same parameters of the previous test, but with  $\nu=10^{-4}$.
In this respect, our method seems very robust and we can reach  very small values of  $\nu$ even for quite
coarse grids. Figure \ref{test2} shows the corresponding results. In this case we see that the solution
is not better than Lipschitz and the support of $Du$ and $m$ are disjoint,
as in the Euclidean case (see \cite[Test 2]{acd}).
\\

\begin{figure}[h!]
 \begin{center}
 \begin{tabular}{cc}
  \footnotesize$\,\,\,\,\,\,\,\min=0.003\,\,\,\,\,\,\,\,\,\,\,\max=1.187$ &
 \footnotesize$\,\,\,\,\,\,\,\min=0\,\,\,\,\,\,\,\,\,\,\,\,\,\,\,\,\,\,\max=37.291$ \\
  \quad\includegraphics[width=.275\textwidth]{mfg-colorbar} &
   \quad\includegraphics[width=.275\textwidth]{mfg-colorbar} \\
\includegraphics[width=.3\textwidth]{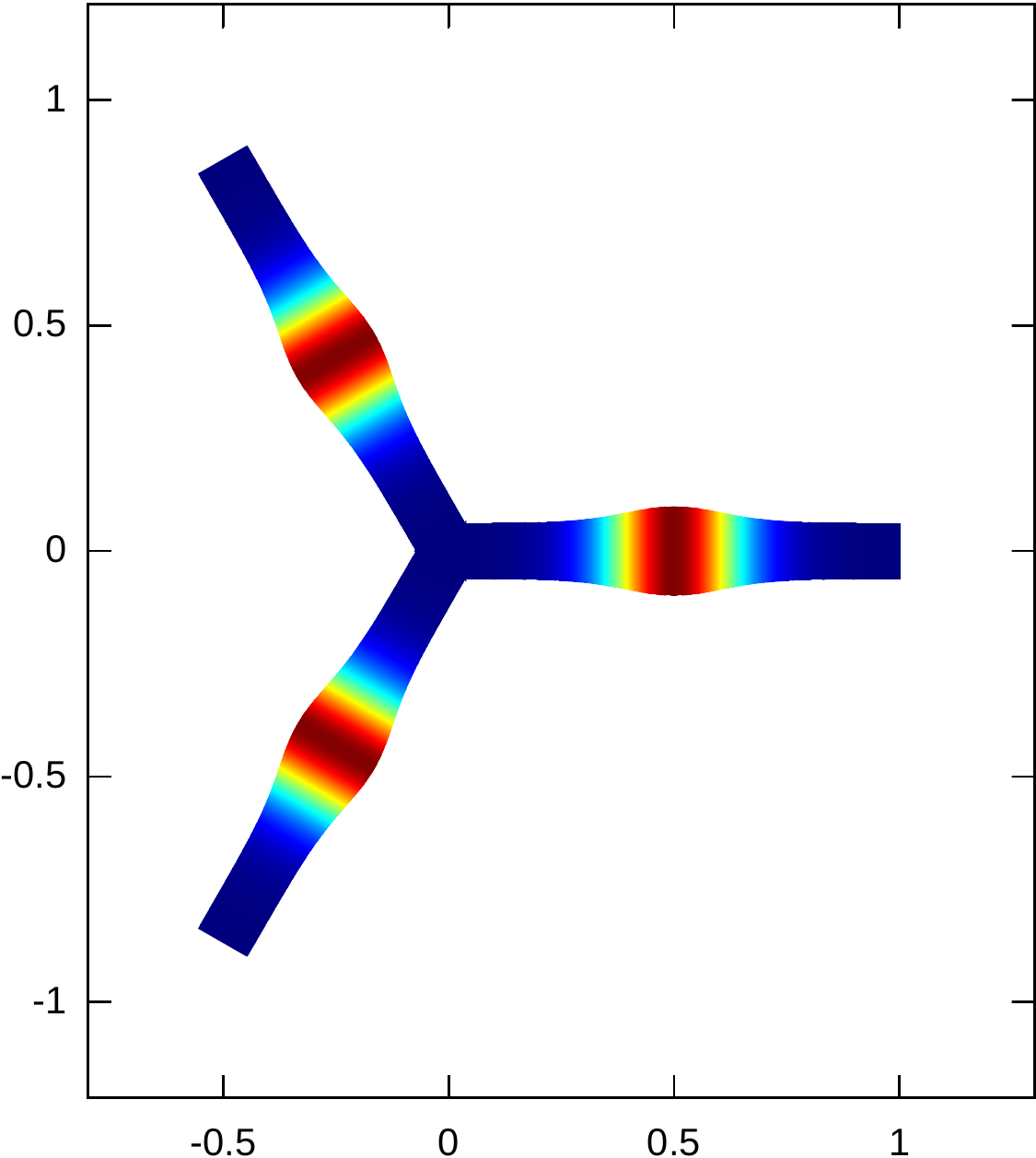} &
\includegraphics[width=.3\textwidth]{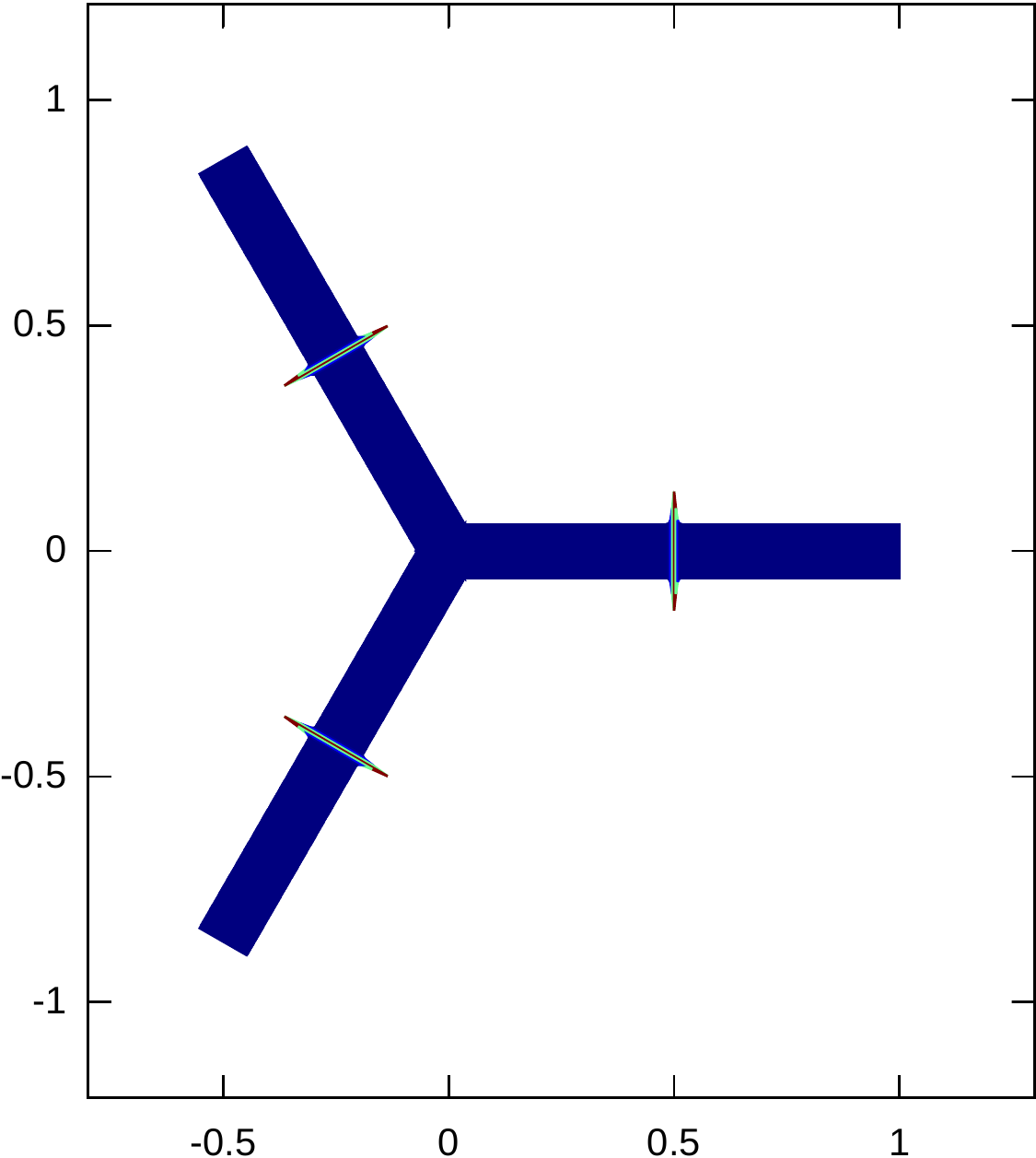} \\
\includegraphics[width=.3\textwidth]{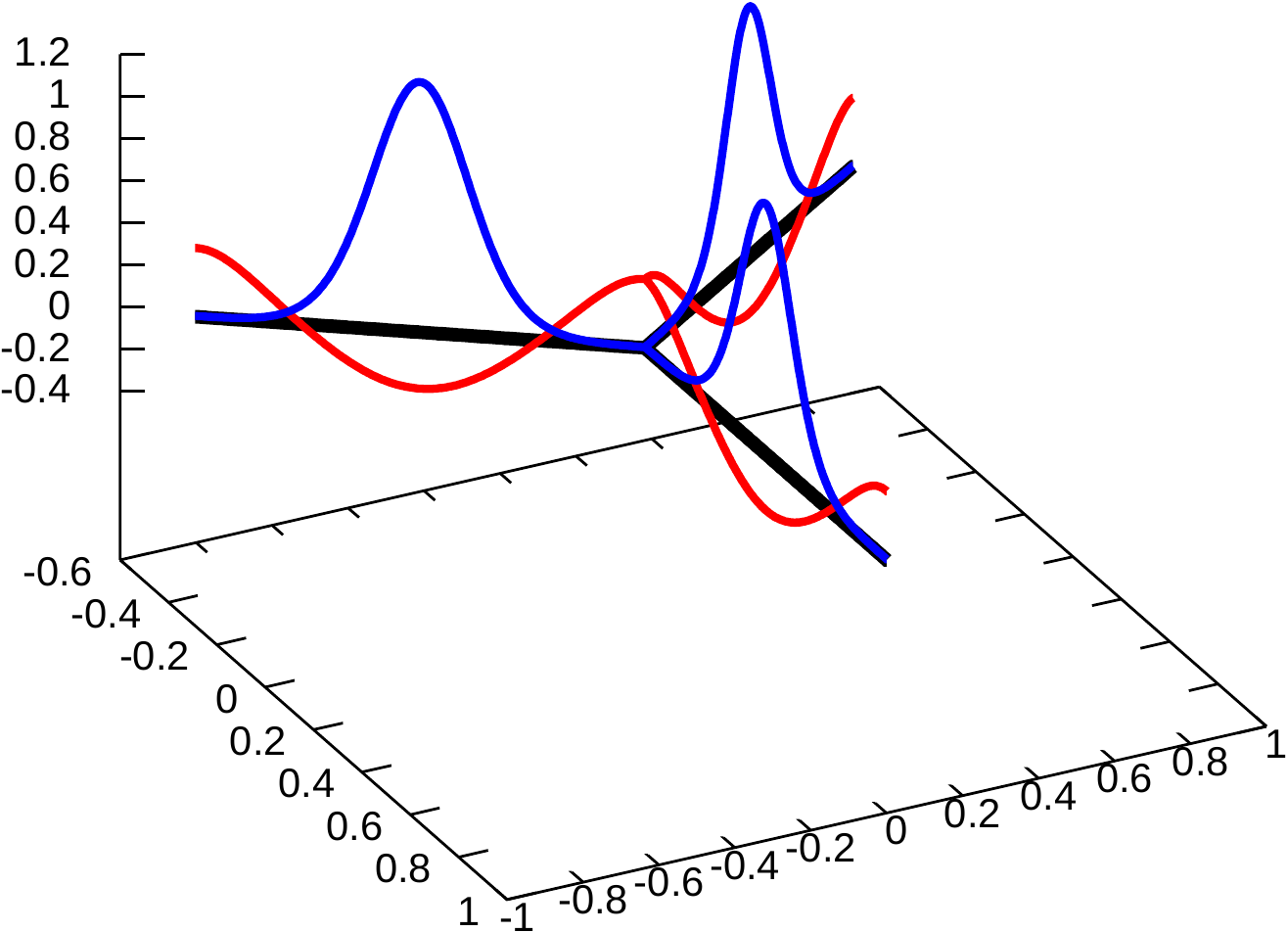} &
\includegraphics[width=.3\textwidth]{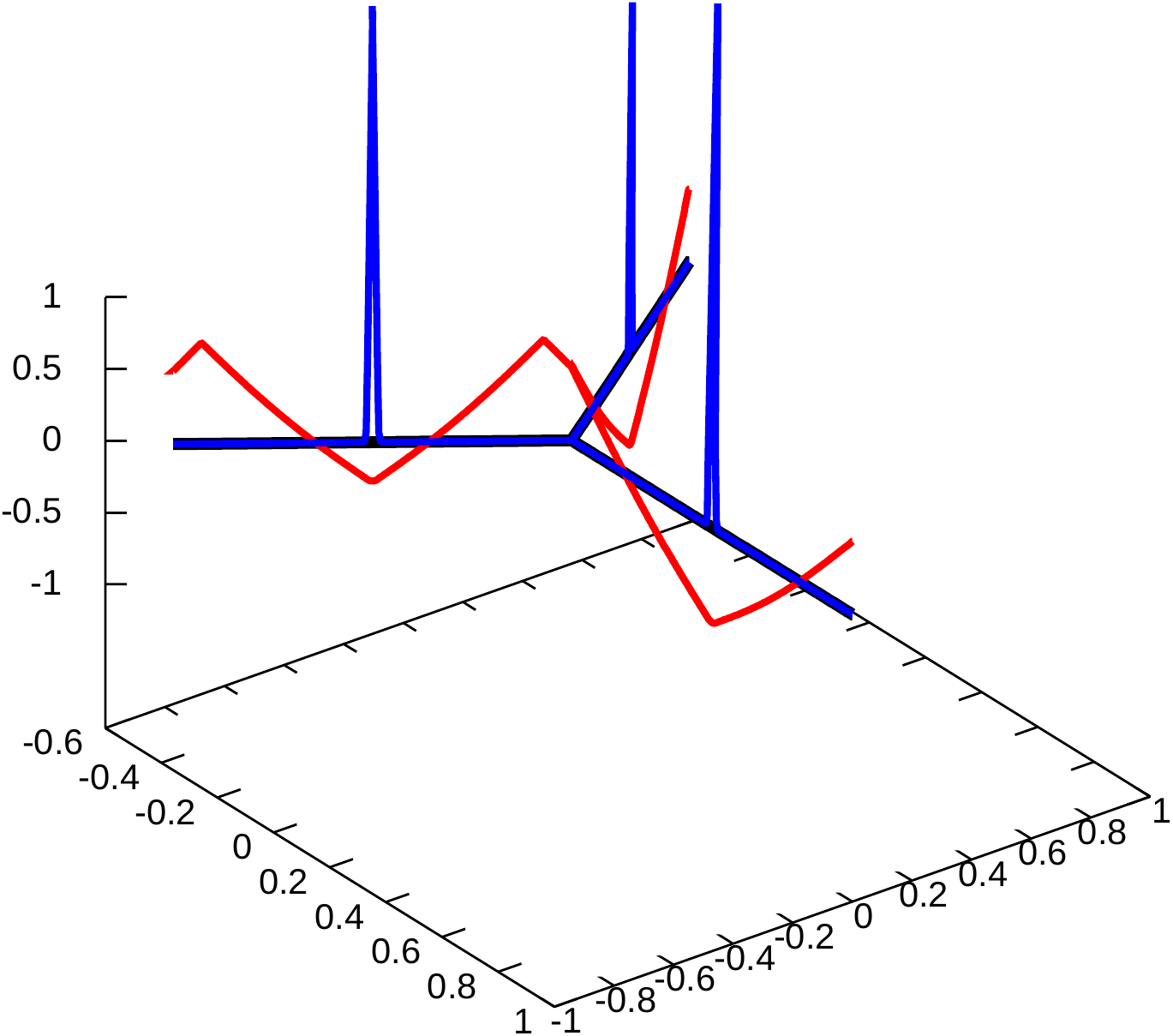} \\
(a)&(b)
\end{tabular}
\end{center}
\caption{the case $V_h[m]=1-\frac{4}{\pi}\arctan(m)$ with (a) $\nu=0.1$ and (b) $\nu=10^{-3}$,
the cost $f$ is active on the whole network.}\label{test3}
\end{figure}

\noindent{\bf Test 3.} We   set $V[m]=1-\frac{4}{\pi}\arctan(m)$, we consider both  $\nu=0.1$ and $\nu=10^{-3}$
and the cost  $f$ active on the whole network, i.e. $(s_0,s_1,s_2)=(1,1,1)$.
Figure \ref{test3} shows the corresponding results. Since $V$ is decreasing, the agents want  the share the
same position and therefore tend  to concentrate around the minima of the value function. Note that for
$\nu$ small, the  regularizing effect of the diffusion is small and  $m$ is close to a sum of  Dirac functions
concentrated at the minima of $u$. In this case assumption \eqref{strictlymon} is not satisfied
and  uniqueness of the solution may fail.
\\

\begin{figure}[h!]
 \begin{center}
 \begin{tabular}{c}
  \footnotesize$\,\,\,\,\,\,\,\min=0\,\,\,\,\,\,\,\,\,\,\,\,\,\,\,\,\,\,\,\,\max=0.103$\\
 \quad\includegraphics[width=.275\textwidth]{mfg-colorbar}\\
 \includegraphics[width=.8\textwidth]{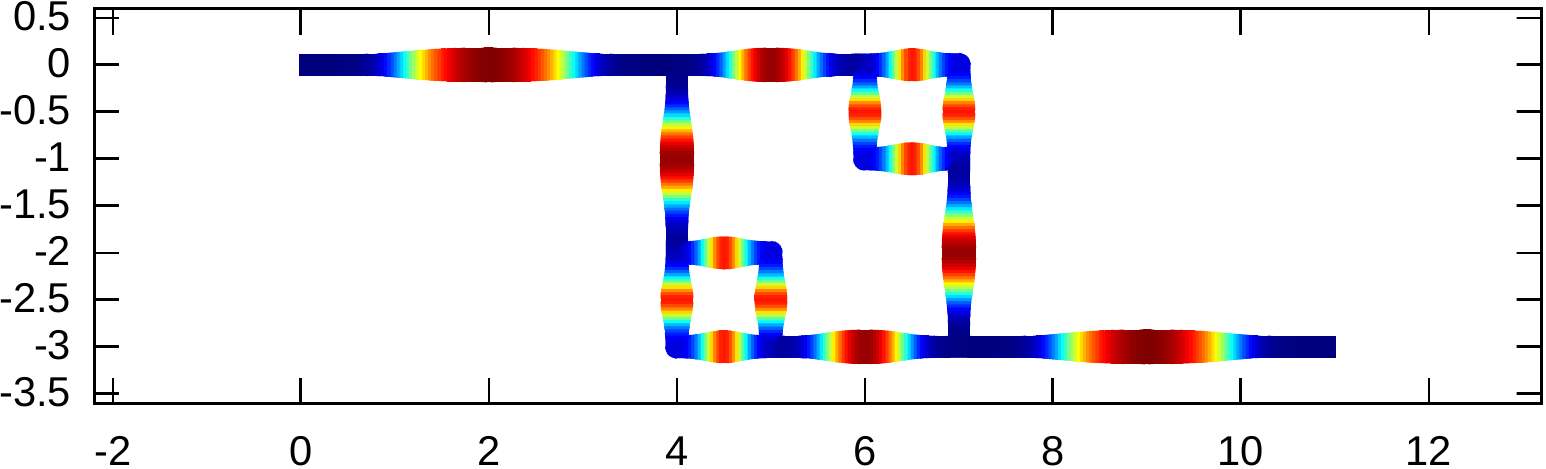} \\(a)\\
 \begin{tabular}{cc}
\includegraphics[width=.4\textwidth]{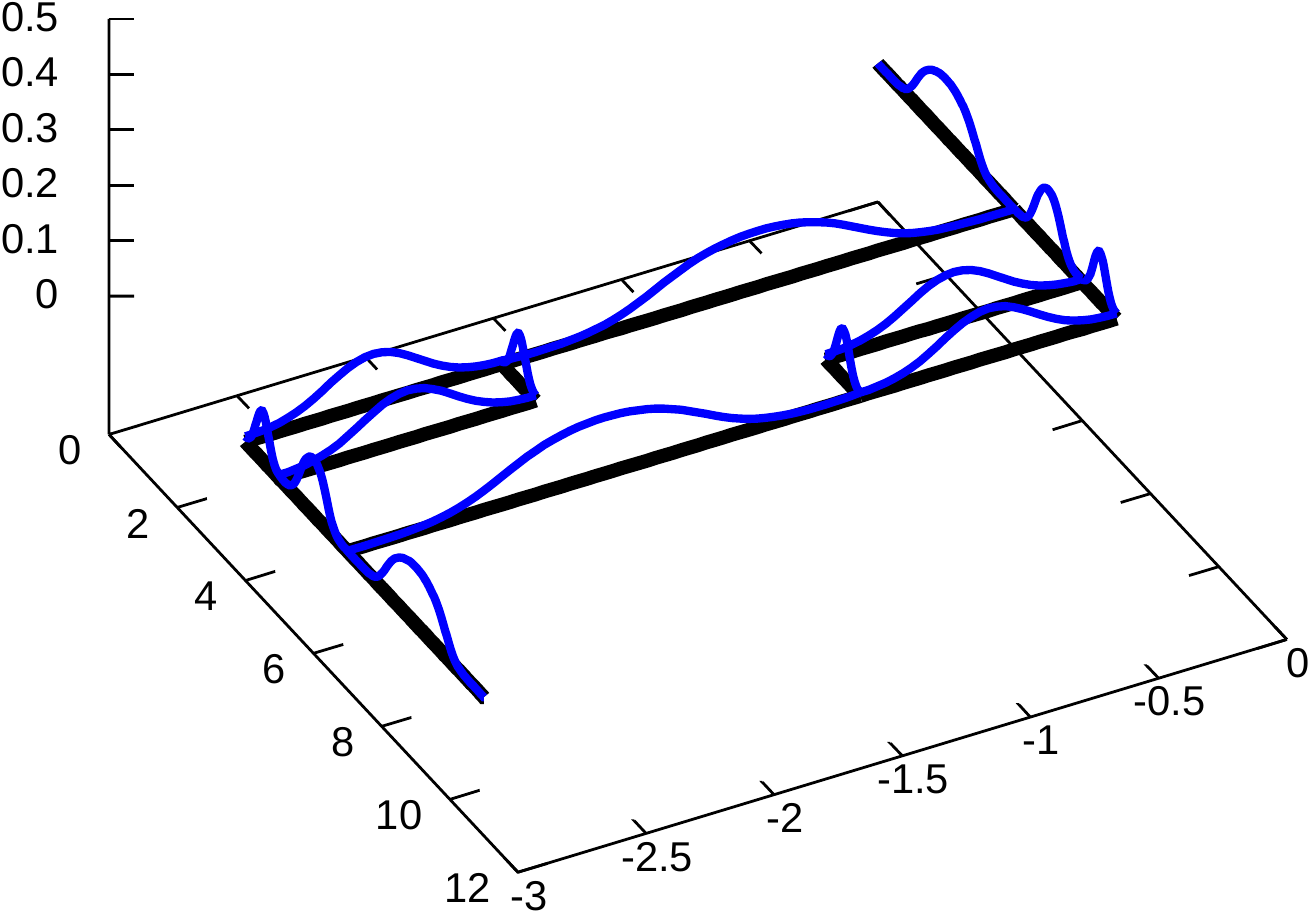} &
\includegraphics[width=.4\textwidth]{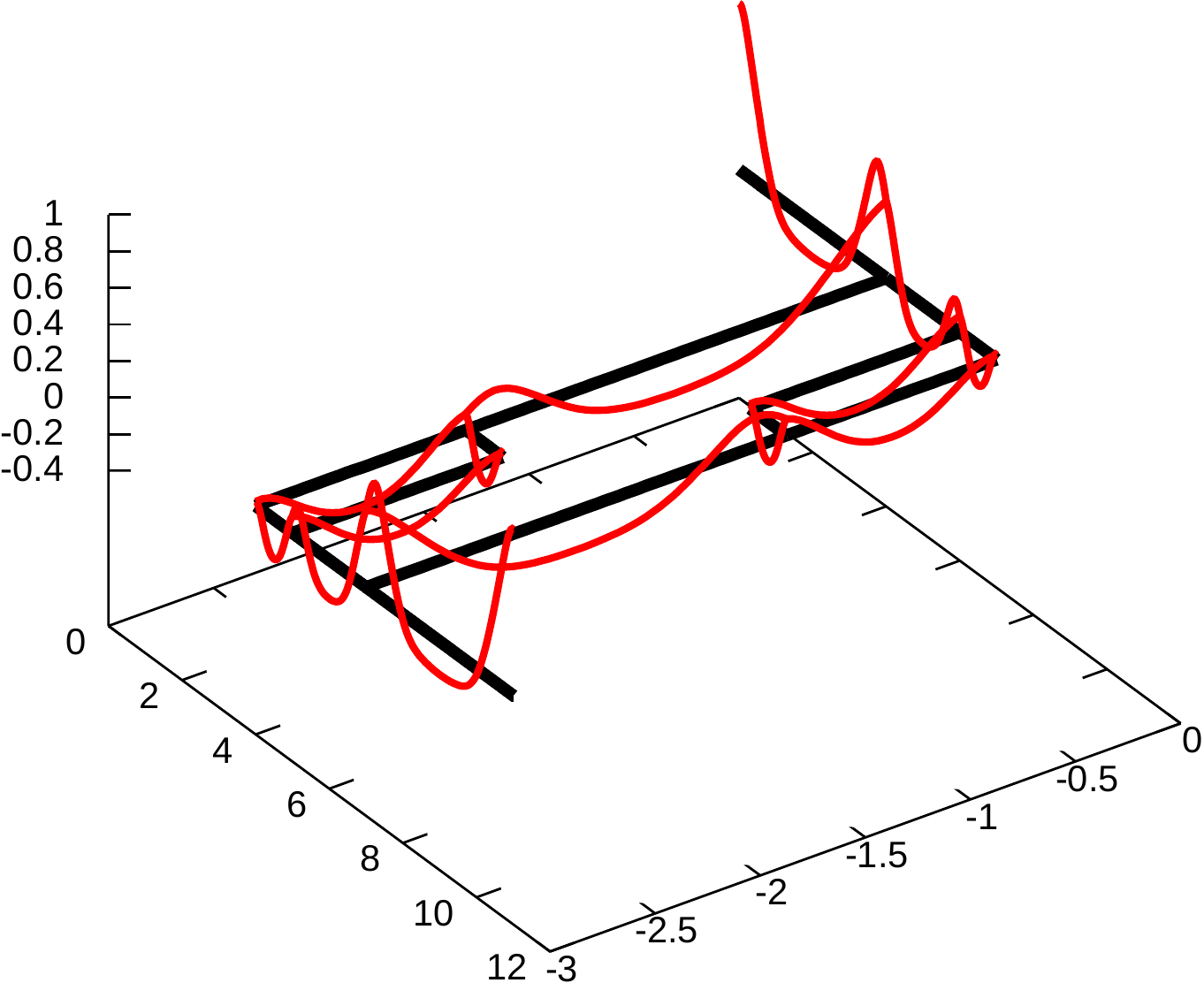} \\
(b)&(c)
\end{tabular}
\end{tabular}
\end{center}
\caption{solution on a self-similar network, (a),(b) the mass $M$, (c) the value function $U$.}\label{test4}
\end{figure}

\noindent{\bf Test 4.} In this experiment, we show that the  method
can efficiently handle the computation on much more complicated structures.
To this end, we consider the periodic network shown in Figure \ref{test4}a. It is a self-similar set, in which
the length of each edge scales with a factor $1/2$ when moving to adjacent edges. Starting from the longest edges,
we stop at the second level of branching and we identify the extremal boundary vertices.
Moreover, we choose the local operator $V_h[m]=m^2$, uniform diffusion coefficient $\nu=0.1$
and the cost  $f$ as before, active on the whole network. In this example  the players are distributed on all the edges with
a scaling factor which depends on the length of the edge.\\

\noindent{\bf Performance and convergence.}
\textbf{Here we present some results showing the convergence and performance of the proposed method, both in terms of accuracy and computational times.
We consider the same setting of Test 1, with the cost $f$ active on all the three edges of the network, i.e. $(s_0,s_1,s_2)=(1,1,1)$.
Moreover, we choose the same number of discretization nodes for each edge, namely $N_j=N$ for $j=0,1,2$ and a variable $N$, so that the space step
is $h=1/N$ on the whole network.
Note that, in the present case, the total number of degrees of freedom ({\em dofs}) of the problem is much more than
$N$. Indeed, we have $N$ nodes for each of the three edges and for both $U$ and $M$, that is {\em dofs}$=6N$.
Since the exact solution is unknown for this problem, we assume as correct the solution computed for $N=2000$, denoted by $(U^{ex},M^{ex},\Lambda^{ex})$.
Then we define the error as
$$
E_h=\|U-U^{ex}\|_1+\|M-M^{ex}\|_1+|\Lambda-\Lambda^{ex}|\,,
$$
where the discrete $1-$norm, for a generic vector $W$ with $3N$ components, is computed as $\|W\|_1=h\sum_{k=1}^{3N}|W_k|$ and the exact solution is projected
on the corresponding grid via linear interpolation. Finally, we define the experimental order of
convergence as {\em Eoc}$(h_1,h_2)=\log(E_{h_1}/E_{h_2})/\log({h_1}/{h_2})$ and we set $\varepsilon=10^{-8}$ for the stopping criterion of the algorithm.}

\begin{figure}[h!]
 \begin{center}
 \begin{tabular}{cc}
\includegraphics[width=.45\textwidth]{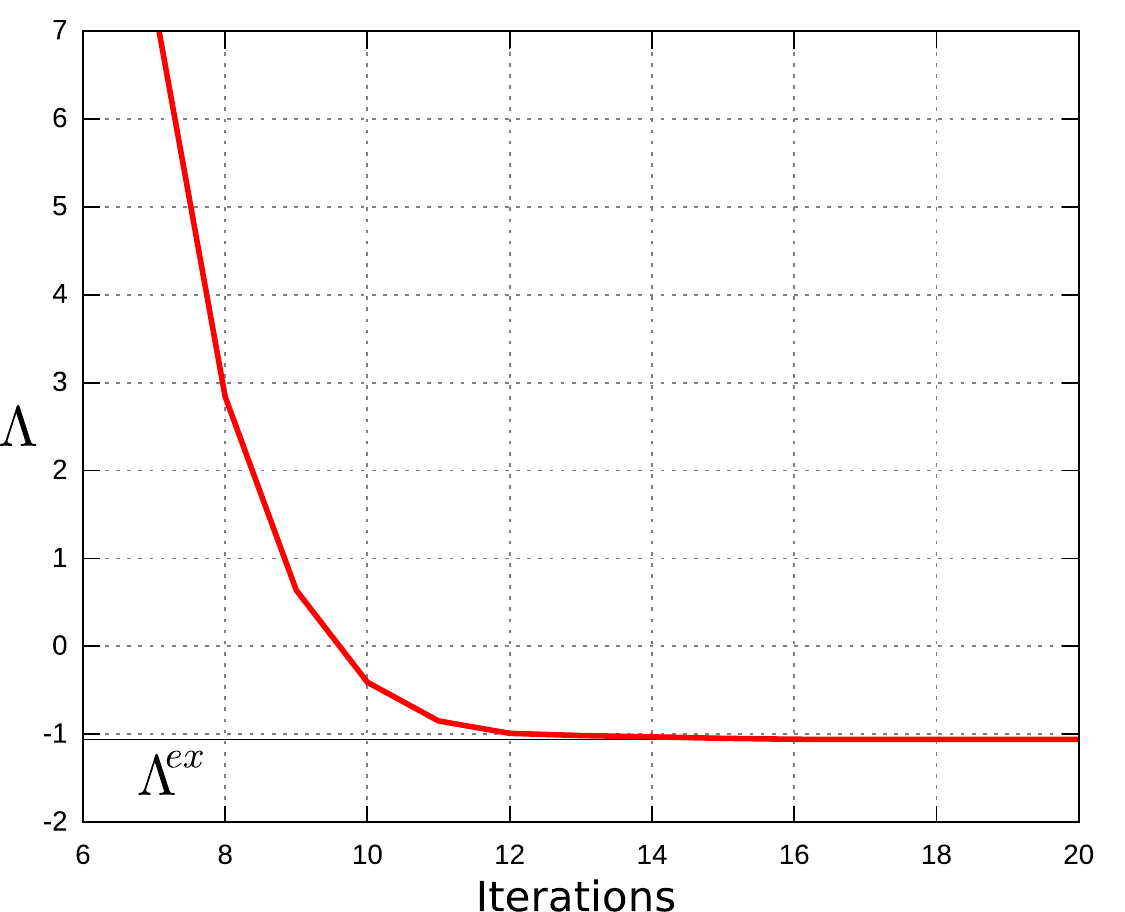} &
\includegraphics[width=.45\textwidth]{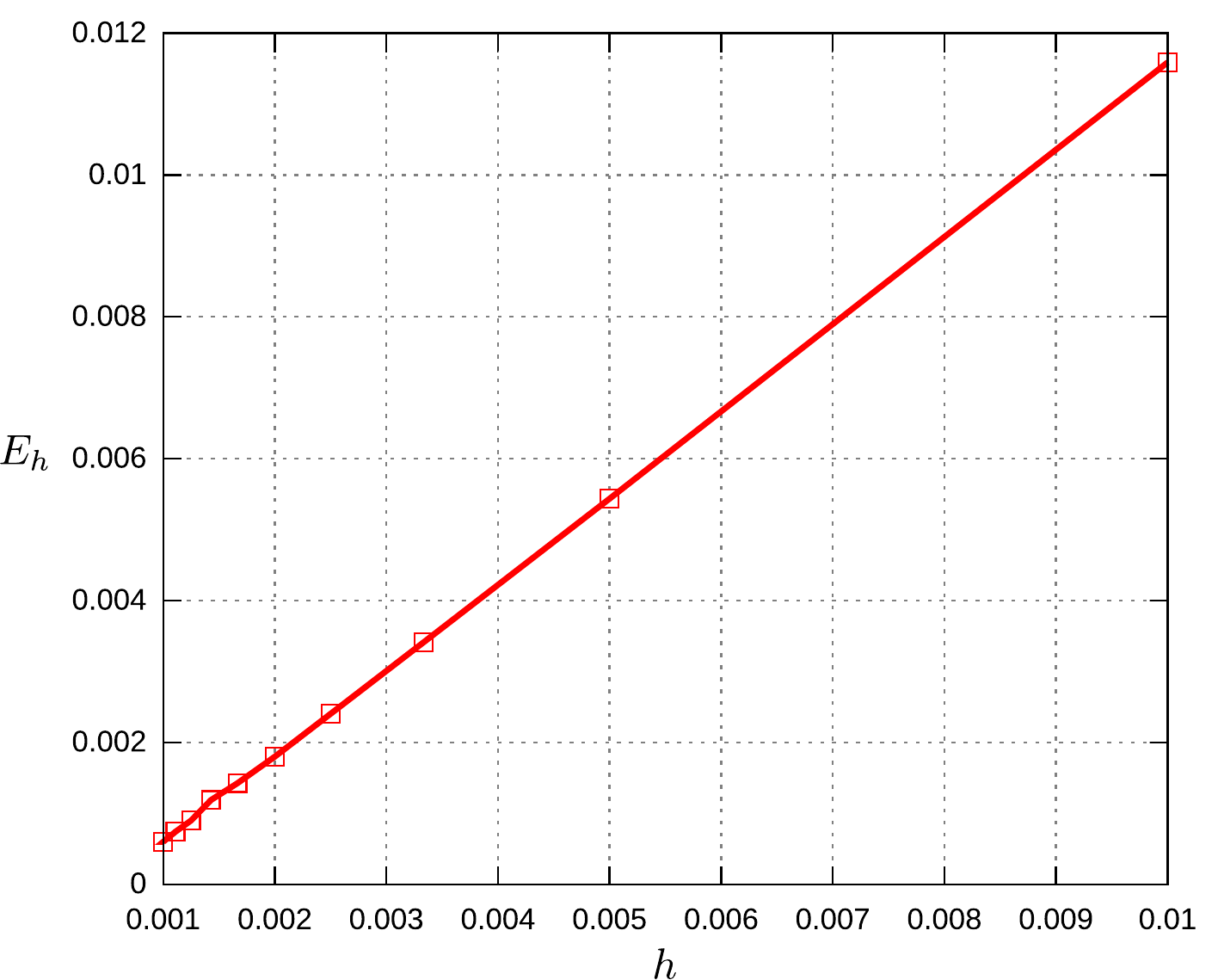} \\
(a)&(b)
\end{tabular}
\end{center}
\caption{convergence, (a) $\Lambda$ vs number of iterations for $N=1000$, (b) error $E_h$ vs space step $h$.}\label{error-test}
\end{figure}

{\bf In Figure \ref{error-test}a we show, for $N=1000$, the behavior of the computed $\Lambda$ as a function of the number of iterations. In this case
$\Lambda^{ex}=-1.058687$, whereas $\Lambda=-1.058876$ is obtained after $20$ iterations with $|\Lambda-\Lambda^{ex}|=0.000189$. 

Similarly, in Figure \ref{error-test}b we plot the error $E_h$ for different space steps $h$, ranging from $10^{-2}$ to $10^{-3}$.
This shows an experimental convergence at least of order $1$. 

Finally, in Table \ref{table-err} we report all the results, including
the error $|\Lambda-\Lambda^{ex}|$ related only to the approximation of the ergodic constant, the {\em Eoc} computed for
successive space steps, the number of iterations and the corresponding computational times. We clearly see that, even for quite coarse grids,
we get a reasonable approximation of the solution with a very low time consumption.}

\begin{table}[!h]
 \centering
  \begin{tabular}{|c|c|c|c|c|c|c|}
    \hline
    N & Dofs & Error $E_h$ & Error $|\Lambda-\Lambda^{ex}|$ & Iterations & Eoc & Cpu time\\
    \hline
    100&600  &0.01159&0.003737& 7&--  & 0.13\\\hline
    200&1200 &0.00544&0.001734& 7&1.09& 0.37\\\hline
    400&2400 &0.00241&0.000762&17&1.17& 4.09\\\hline
    800&4800 &0.00091&0.000284&16&1.40&19.57\\\hline
    1000&6000&0.00059&0.000189&20&1.94&47.94\\\hline
 \end{tabular}
  \caption{Performance of the proposed method.}\label{table-err}
\end{table}


\end{document}